\documentclass[review]{elsarticle}
\usepackage[margin=1.5cm]{geometry}
\usepackage{lineno,hyperref}
\usepackage{amssymb}
\usepackage{amsmath}
\usepackage{tabularx}
\usepackage[english]{babel}
\usepackage{graphicx}
\usepackage{color}
\usepackage{epstopdf}
\usepackage{multirow}
\usepackage{subfig} 
\usepackage{bm}
\usepackage{algorithm}
\usepackage{threeparttable}
\usepackage{bigints}

\modulolinenumbers[5]

\newdefinition{rmk}{Remark}
\newdefinition{define}{Definition} 
\newproof{pf}{Proof}

\journal{Elsevier}

\begin{document}

\begin{frontmatter}

\title{A consistent analytical formulation for volume-estimation of geometries enclosed by implicitly defined surfaces}

\author{Shucheng Pan} \ead{shucheng.pan@tum.de}
\author{Xiangyu Hu} \ead{xiangyu.hu@tum.de}
\author{Nikolaus. A. Adams} \ead{nikolaus.adams@tum.de}
\address{Lehrstuhl f\"{u}r Aerodynamik und Str\"{o}mungsmechanik, Technische Universit\"{a}t
M\"{u}nchen, 85748 Garching, Germany}
\begin{abstract}
We have derived an analytical formulation for estimating the volume of geometries enclosed by implicitly defined surfaces. 
The novelty of this work is due to two aspects. 
First we provide 
a general analytical formulation for all two-dimensional cases,
and for elementary three three-dimensional cases 
by which the volume of general three-dimensional cases can be computed. 
Second, our method addresses the inconsistency issue due to mesh refinement.
It is demonstrated by several two-dimensional and three-dimensional cases that this analytical formulation exhibits 2nd-order accuracy.
\end{abstract}

\end{frontmatter}


\section{\label{sec:intro}Introduction}
Accurate and efficient volume estimation is important for investigating complex geometries encountered 
in a variety of scientific and engineering problems. 
In astronomy, the estimation of an asteroid volume is helpful to determine 
the bulk density for assessing the material composition. 
Except for some large spheroidal asteroids whose volumes can be estimated by their mean diameters, 
most asteroids have irregular shapes which require more complicated volume estimation methods. 
Usually, three-dimensional (3D) reconstructed shape models are employed to approximate 
such irregular geometries \cite{fujiwara2006rubble, besse20093}. 
In computational fluid dynamics, 
a volume-fraction estimation is required for the transport of multiple fluids 
by the volume-of-fluid (VOF) method or the level-set method \cite{hu2006conservative}. 
The volume fraction of a specific fluid in a uniform mesh cell 
can be computed by the Heaviside \cite{so2011anti, Sussman2000voflevelset} or
color function \cite{gueyffier1999volume}, 
or combination of either one with the interface normal direction \cite{hu2010multi}. The interface reconstruction step of the VOF method may also require a volume-estimation method to determine the interface location \cite{scardovelli2000analytical, diot2014interface, diot2016interface}.
Another example is medical engineering, where in order to detect the development of lung hypoplasia 
and investigate the correlation between lung growth and fetal presentation, 
accurate measurement of fetal lung volume is a pressing need \cite{rypens2001fetal}. 
Tumor volume estimation is also widely used in many cancer treatments, 
such as prostate cancer \cite{coakley2002prostate}, brain tumors \cite{joe1999brain}
and pelvic neoplasms \cite{mayr1996usefulness}. 
Most of such volume estimations are based on magnetic-resonance (MR) imaging. 
A three-dimensional reconstruction method is considered to be much more reliable 
and accurate than other more traditional approximation methods, 
such as ellipsoid and disc-summation approximations, 
and automated image segmentation \cite{Durso2014, Gong2012, Cheong2007, Jones1983}. 
For example, volume estimation of the articular cartilage is performed 
by the Marching Cube algorithm \cite{Lorensen1987marchingcube, Newman2006} applied on the three-dimensional reconstruction \cite{nystroem2002area}. 
The disadvantage of this method is that one first needs to reconstruct the iso-surface 
before volume estimation. However, there are too many different cases 
($16$ for 2D and $256$ for 3D) to deal with. The volume estimation method is consistent during mesh refinement in the sense that the volume fractions are conservative with respect to the coarsest level. It is possible to obtain a unique volume approximation during mesh refinement by a consistent volume estimation method. Any numerical method introduces inconsistencies during mesh refinement, irrespectively of its order of accuracy, as the truncation error depends on mesh resolution. Although not critical for many applications, consistency is desirable, and it is important for numerical simulations of multiphase flows \cite{sussman1999adaptive}.

The objective of the present paper is to develop a consistent volume estimation method for calculating the volumes of geometries enclosed by implicitly defined surfaces. To achieve this property, we derive analytical formulations for generic 2D and 3D geometries, which avoid truncation errors, unlike numerical methods. The formulations are then subsumed to derive a general formulation that can be applied to more complex cases. This method is made consistent during mesh refinement by employing bilinear/trilinear interpolation. The paper is organized as follows. Sec. \ref{sec:Analy} describes the derivation of 2D and 3D analytical volume formulation for arbitrary geometries. We also discuss how to preserve consistency during mesh refinement. Sec. \ref{sec:app} is dedicated to assessing the capability of the present method to calculate the volumes for various cases, followed by conclusions in Sec. \ref{sec:con}.

\section{\label{sec:Analy}Analytical Volume Fraction}

\subsection{2D formulation}

Before presenting the 3D formulation, we first consider the 2D case for its simplicity and give a comprehensive description of the underlying concept. Suppose that we have volumetric data, say discrete data of a level-set field $\phi$, which implicitly defines the interface in the domain $\Omega$, the interface in any cell can be expressed by the zero value of a reconstructed field $\phi(x,y)$ which is obtained by bilinear interpolation. The interface is hereby represented by piecewise parabolas for two dimensions (2D), i.e., the interface inside the cell $[i,i+1]\times[j,j+1]$ is defined as
\begin{equation}\label{eqvar2dphi}
{{\Gamma }}=\left\{ \left( x,y \right)\in \Omega :\phi(x,y)=\beta_0+\beta_1x+\beta_2y+\beta_3xy=0 \right\}, 
\end{equation}
where the coefficients $\beta_0$, $\beta_1$, $\beta_2$ and $\beta_3$ are determined by the interpolation conditions,
\begin{eqnarray}\label{eqvar2d}
\beta_0=\phi_{00}, \quad  \beta_1=\phi_{10}-\phi_{00},\quad  \beta_2=\phi_{01}-\phi_{00}, \quad \beta_3=\phi_{00}+\phi_{11}-\phi_{01}-\phi_{10},
\end{eqnarray}
with $\phi_{00}$, $\phi_{10}$, $\phi_{01}$ and $\phi_{11}$ defined in Fig. \ref{sketch}(a). The volume of this cell can be easily calculated by integrating the reconstructed field in $x$ direction or $y$ direction. In order to simplify the integration, we reformulate Eq. (\ref{eqvar2dphi}) with a coordinate mapping,
\begin{equation}\label{eqvar2d_new}
\zeta=\frac{-a \eta-b}{c \eta+d}
\end{equation}
whose symbols need careful definition to simplify the calculations. The definition is based on the intersection points at the edges of a mesh cell. Specifically in 2D, we first detect the grid point
\begin{equation}\label{origin}
\mathbf{x}_0 = \underset{\mathbf{x} \in \mathbf{X}}{\text{arg max }} \# \left( \mathbf{x}' \in \mathbf{X} | \phi(\mathbf{x}) \phi(\mathbf{x}') < 0 \wedge  \|\mathbf{x}' - \mathbf{x} \|_2 = 1\right)
\end{equation}
and set is as the origin of a local Cartesian coordinate system. $\mathbf{X}$ is the set of all $4$ vertices ($8$ in 3D) of the current cell. Then the coordinate transformation is determined by $\eta = \underset{\gamma \in \{ x, y  \}}{\text{arg max }} l_{\gamma}$ and $\zeta = \left\lbrace  x,y \right\rbrace  \backslash \left\lbrace \eta\right\rbrace$, where $l_{\gamma}$ is the distance from the vertex $(i,j)$ to the intersection point on a neighbour edge in $x$-direction or $y$-direction,
\begin{eqnarray}\label{eq1ength2d}
{{l}_{x}} = -\frac{{{\beta }_{2}}{{i}_{y}}+{{\beta }_{0}}}{{{\beta }_{3}}{{i}_{y}}+{{\beta }_{1}}}  \quad  \mbox{if} \quad {{\phi }_{{{i}_{x}}{{i}_{y}}}}{{\phi }_{|1-{{i}_{x}}|{{i}_{y}}}}<0 \quad \mbox{else} \quad 1.0, \quad
{{l}_{y}} = -\frac{{{\beta }_{1}}{{i}_{x}}+{{\beta }_{0}}}{{{\beta }_{3}}{{i}_{x}}+{{\beta }_{2}}}  \quad  \mbox{if} \quad {{\phi }_{{{i}_{x}}{{i}_{y}}}}{{\phi }_{{{i}_{x}}|1-{{i}_{y}}|}}<0 \quad \mbox{else} \quad 1.0,
\end{eqnarray}
where $i_x$ and $i_y$ are the coordinates of the origin in the $x$ and $y$ directions, respectively, see Fig. \ref{sketch}(a). Once the origin and the local transformed coordinates ($\eta, \zeta$) are defined, the coefficients of Eq. (\ref{eqvar2d_new}), $a$, $b$, $c$ and $d$, are determined according to the relation between the $(x,y)$ and $(\eta,\zeta)$ coordinate systems. Then the cut volume of cell $[i,i+1]\times[j,j+1]$ can be obtained by integration of $\zeta(\eta)$ as
\begin{eqnarray}\label{volume2d}
\alpha  = \left. (1-2i_{\varsigma})\left[ \frac{(ad-bc)\log|c\eta+d|}{c^2}-\frac{a\eta}{c}\right] \right|_{\eta_0}^{\eta_1}
+i_{\varsigma}({{\eta }_{1}} - {{\eta }_{0}}).
\end{eqnarray}
The integration range is determined by  $\eta_{0} = \mbox{min}(l_{\eta}, i_{\eta})$ and $\eta_{1} = \mbox{max}(l_{\eta}, i_{\eta})$, where ${l}_{\eta}$ is the length from the origin of the coordinate system to the intersection point in $\eta$ direction.

\subsection{3D formulation}
\begin{figure}[p]
\begin{center}
\includegraphics[width=0.8\textwidth]{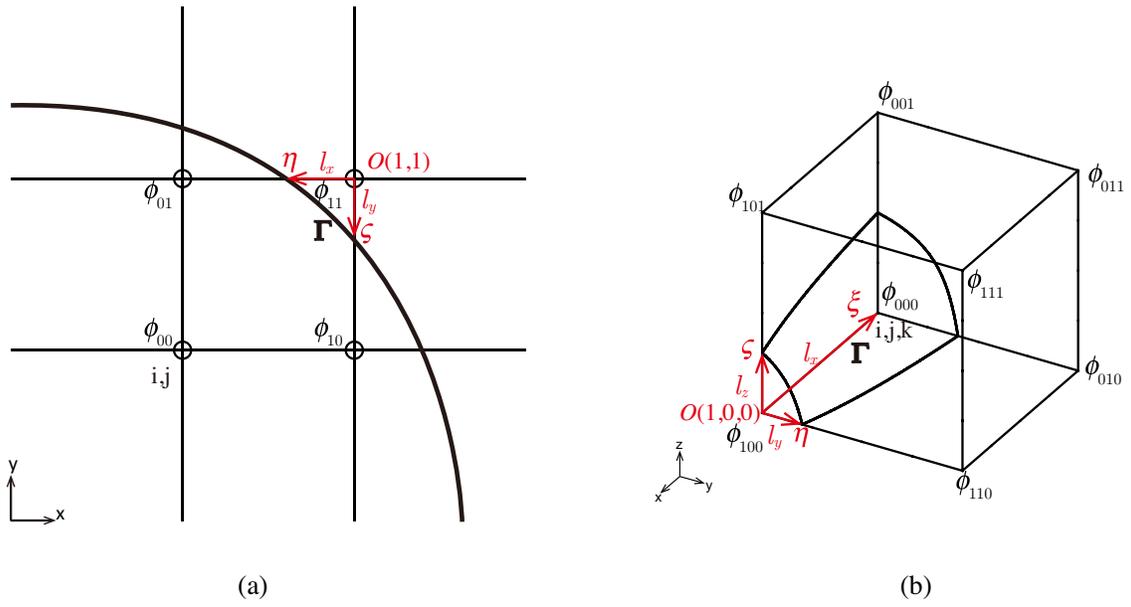}
\caption{(a) 2D and (b) 3D volume estimation for the cell $[i,i+1]\times[j,j+1]$ and $[i,i+1]\times[j,j+1]\times[k,k+1]$, respectively.}
\label{sketch}
\end{center}
\end{figure}

The 2D formulation can be easily extended to three dimensions. Again, the interface inside a cell, say $[i,i+1]\times[j,j+1]\times[k,k+1]$ in Fig. \ref{sketch}(b), can be represented by trilinear interpolation,
\begin{eqnarray}\label{eq3d}
{{\Gamma}} = \left\{\left( x,y,z \right)\in \Omega :{{\phi }}(x,y,z) \right. = \left.{{\beta }_{0}}+{{\beta }_{1}}x+{{\beta }_{2}}y+{{\beta }_{3}}z+{{\beta }_{4}}xy+{{\beta }_{5}}yz+{{\beta }_{6}}xz+{{\beta }_{7}}xyz=0 \right\},
\end{eqnarray}
where the coefficients are uniquely determined by the interpolation conditions, see Fig. \ref{sketch}(b),
\begin{eqnarray}\label{eqvar}
&&\beta_0=\phi_{000}, \qquad  \beta_1=\phi_{100}-\phi_{000},\qquad  \beta_2=\phi_{010}-\phi_{000}, \qquad \beta_3=\phi_{001}-\phi_{000}, \nonumber \\
&&\beta_4=\phi_{110}-\phi_{100}-\phi_{010}+\phi_{000}, \qquad \beta_5=\phi_{011}-\phi_{010}-\phi_{001}+\phi_{000}, \nonumber \\
&&\beta_6=\phi_{101}-\phi_{100}-\phi_{001}+\phi_{000}, \qquad \beta_7=\phi_{111}-\phi_{110}-\phi_{101}+\phi_{100}-\phi_{011}+\phi_{010}+\phi_{001}-\phi_{000}.
\end{eqnarray}
Upon identification of the origin ($O$) by Eq. (\ref{origin}), the local coordinates are defined as $\xi = \underset{\gamma \in \{ x, y, z \}}{\text{arg max }} l_{\gamma}$, $\eta = \underset{\gamma \in \{ x, y, z \} \setminus \{\xi\}}{\text{arg max }} l_{\gamma}$ and $\zeta = \{x,y,z\} \setminus \{\xi, \eta\}$, where the distance from the grid $(i,j,k)$ to the intersection points, $l_{\gamma}$, are calculated by
\begin{eqnarray}\label{eq1ength3d}
{{l}_{x}} &=& \frac{-\left[ \left( {{\beta }_{5}}{{i}_{y}}+{{\beta }_{3}} \right){{i}_{z}}+{{\beta }_{2}}{{i}_{y}}+{{\beta }_{0}} \right]}{\left[ \left( {{\beta }_{7}}{{i}_{y}}+{{\beta }_{6}} \right){{i}_{z}}+{{\beta }_{4}}{{i}_{y}}+{{\beta }_{1}} \right]}  \quad  \mbox{if} \quad {{\phi }_{{{i}_{x}}{{i}_{y}}{{i}_{z}}}}{{\phi }_{|1-{{i}_{x}}|{{i}_{y}}{{i}_{z}}}}<0 \quad \mbox{else} \quad 1.0 \nonumber \\
{{l}_{y}} &=& \frac{-\left[ \left( {{\beta }_{6}}{{i}_{x}}+{{\beta }_{3}} \right){{i}_{z}}+{{\beta }_{1}}{{i}_{x}}+{{\beta }_{0}} \right]}{\left[ \left( {{\beta }_{7}}{{i}_{x}}+{{\beta }_{5}} \right){{i}_{z}}+{{\beta }_{4}}{{i}_{x}}+{{\beta }_{2}} \right]}  \quad  \mbox{if} \quad {{\phi }_{{{i}_{x}}{{i}_{y}}{{i}_{z}}}}{{\phi }_{{{i}_{x}}|1-{{i}_{y}}|{{i}_{z}}}}<0 \quad \mbox{else} \quad 1.0 \nonumber \\
{{l}_{z}} &=& \frac{-\left[ \left( {{\beta }_{4}}{{i}_{x}}+{{\beta }_{2}} \right){{i}_{y}}+{{\beta }_{1}}{{i}_{x}}+{{\beta }_{0}} \right]}{\left[ \left( {{\beta }_{7}}{{i}_{x}}+{{\beta }_{5}} \right){{i}_{y}}+{{\beta }_{6}}{{i}_{x}}+{{\beta }_{3}} \right]}  \quad  \mbox{if} \quad {{\phi }_{{{i}_{x}}{{i}_{y}}{{i}_{z}}}}{{\phi }_{{{i}_{x}}{{i}_{y}}|1-{{i}_{z}}|}}<0 \quad \mbox{else} \quad 1.0.
\end{eqnarray}
The coordinates of the origin in the $x$, $y$ and $z$ directions are $i_x$, $i_y$ and $i_z$, respectively, as shown in Fig. \ref{sketch}(b). Eq. (\ref{eq3d}) is rewritten as
\begin{equation}
\zeta(\xi, \eta) =\frac{-\xi(a \eta+b)-c \eta-d}{\xi(e \eta+f)+g \eta+h},
\end{equation}
with the coefficients, $a, b, \dots, f, g$, determined by the relation between the $(x,y,z)$ and $(\xi,\eta,\zeta)$ coordinate systems. Thus the volume of each cell $\alpha$ can be calculated by the double integral of $\zeta(\xi,\eta)$,
\begin{eqnarray}\label{volume3d}
\alpha=(1-2i_{\varsigma})\int_{{{\xi }_{0}}}^{{{\xi }_{1}}}{\left[ \int_{{{\eta }_{0}}}^{{{\eta }_{1}}}{\zeta(\xi ,\eta )d\eta } \right]d}\xi +i_{\varsigma}\int_{{{\xi }_{0}}}^{{{\xi }_{1}}}q(\xi)d\xi = (1-2i_{\varsigma})\left[ \mathbf{F}\left( \xi_{1}  \right)- \mathbf{F}\left( \xi_{0}  \right)\right]
+i_{\varsigma}\left[ \mathbf{G}\left( \xi_{1}  \right)- \mathbf{G}\left( \xi_{0}  \right)\right]
\end{eqnarray}
with
\begin{equation}
q(\xi)=\frac{-(fi_{\varsigma}+b)\xi-hi_{\varsigma}-d}{(ei_{\varsigma}+a)\xi+gi_{\varsigma}+c}.
\end{equation}
The integration ranges are $\xi_{0} = \mbox{min}(l_{\xi}, i_{\xi})$ and $\xi_{1} = \mbox{max}(l_{\xi}, i_{\xi})$ in $\xi$-direction, and
\begin{eqnarray}
\eta_{0} =-\frac{\left( f\,i_{\zeta}+b\right) \,\xi+h\,i_{\zeta}+d}{\left( e\,i_{\zeta}+a\right) \,\xi+g\,i_{\zeta}+c} \quad \mbox{if} \quad i_{\eta}=1 \quad \mbox{else} \quad 0.0, \quad
\eta_{1} = -\frac{\left( f\,i_{\zeta}+b\right) \,\xi+h\,i_{\zeta}+d}{\left( e\,i_{\zeta}+a\right) \,\xi+g\,i_{\zeta}+c} \quad \mbox{if} \quad i_{\eta}=0 \quad \mbox{else} \quad 1.0 ,
\end{eqnarray}
in $\eta$-direction.

\begin{figure}[p]
\begin{center}
\includegraphics[width=0.8\textwidth]{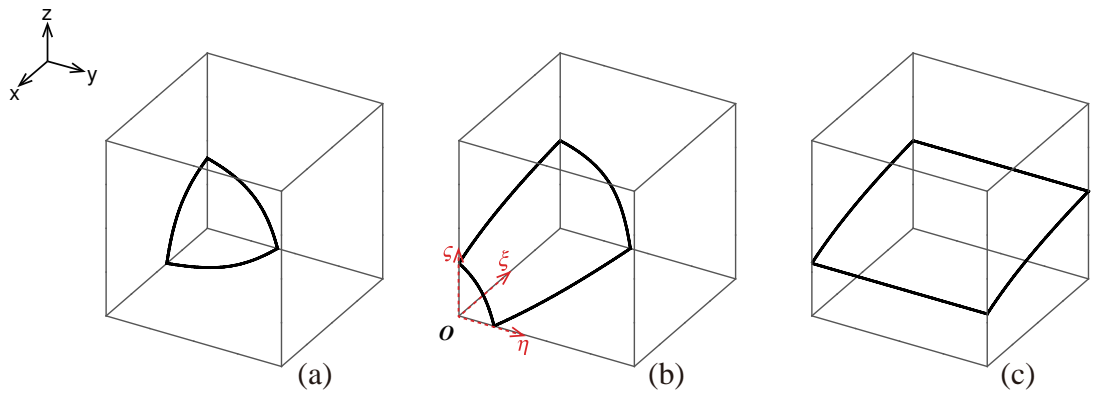}
\caption{Three elementary cases for three dimensional volume estimations.}
\label{3d_basic}
\end{center}
\end{figure}

For all cells containing a resolved interface, Eq. (\ref{volume3d}) for the three cases sketched in Fig. \ref{3d_basic} can be expressed by elementary integrals. First we subsume those elementary cases by a single formulation. Other cases can be expressed as combination of these elementary cases, and their volume estimation will be discussed later. 

First we consider $i_{\zeta} = 0$ such that $\mathbf{F}(\xi)$ in Eq.(\ref{volume3d}) becomes
\begin{eqnarray}
\mathbf{F}(\xi)=\int{\left[ -\frac{\left( {{t}_{2}}+\xi{{t}_{1}}+{{\xi}^{2}}{{t}_{0}} \right)}{{{(e\xi+g)}^{2}}}
\left[\log \left( f\xi+h \right)+\log \left( a\xi+c \right) -\log \left( {{t}_{2}}+\xi{{t}_{1}}+{{\xi}^{2}}{{t}_{0}} \right) \right]+\frac{b\xi+d}{e\xi+g} \right]d\xi}
\end{eqnarray}
with $t_0=af-be$, $t_1=ah-bg+cf-de$, and $t_2=ch-dg$, while $\mathbf{G}(\xi)$ is
\begin{eqnarray}\label{G}
{\mathbf{G}}(\xi)=-\frac{1}{({{i}_{\zeta }}e+a)^2}\left[ \left( eh-fg \right)i_{\zeta }^{2}+\left( ah-bg-cf+de \right){{i}_{\zeta }} + ad-bc \right]\log \left(  {{i}_{\zeta }}{{\xi}_{0}}+{{\xi}_{4}}  \right)-\frac{\left( {{i}_{\zeta }}f+b \right)\xi}{{{i}_{\zeta }}e+a}
\end{eqnarray}
We define $\mathbf{F}(\xi) = \mathbf{F}_1(\xi)+\mathbf{F}_2(\xi)$, where $\mathbf{F}_1$ and $\mathbf{F}_2$ are
\begin{eqnarray}
\mathbf{F}_1(\xi)&=&\log({{t}_{2}}+\xi{{t}_{1}}+{{\xi}^{2}}{{t}_{0}})\left[-\frac{{{e}^{2}}{{t}_{2}}-eg{{t}_{1}}+(-{{e}^{2}}{{\xi}^{2}}-eg\xi+{{g}^{2}}){{t}_{0}}}{{{e}^{4}}\xi+{{e}^{3}}g} \right. + \left.\frac{\log (e\xi+g)(e{{t}_{1}}-2g{{t}_{0}})}{{{e}^{3}}}\right]
\nonumber \\
&-&\int\left[-\frac{{{e}^{2}}{{t}_{2}}-eg{{t}_{1}}
+(-{{e}^{2}}{{\xi}^{2}}-eg\xi+{{g}^{2}}){{t}_{0}}}{{{e}^{4}}\xi+{{e}^{3}}g} \right. + \left.\frac{\log (e\xi+g)(e{{t}_{1}}-2g{{t}_{0}})}{{{e}^{3}}}\right]
d({{t}_{2}}+\xi{{t}_{1}}+{{\xi}^{2}}{{t}_{0}}),
\end{eqnarray}
and
\begin{eqnarray}
\mathbf{F}_2(\xi)&=&\int\left[ -\frac{\log \left( f\xi+h \right)\left( {{t}_{2}}+\xi{{t}_{1}}+{{\xi}^{2}}{{t}_{0}} \right)}{{{(e\xi+g)}^{2}}}-\frac{\log \left( a\xi+c \right)\left( {{t}_{2}}+\xi{{t}_{1}}+{{\xi}^{2}}{{t}_{0}} \right)}{{{(e\xi+g)}^{2}}} \right. + \left. \frac{b\xi+d}{e\xi+g} \right]d\xi,
\end{eqnarray}
respectively. Substituting the variables, ${{\Xi}_{0}}=e\xi+g$, ${{\Xi}_{1}}={{t}_{2}}+\xi{{t}_{1}}+{{\xi}^{2}}{{t}_{0}}$, ${{t}_{3}}={{e}^{2}}{{t}_{2}}-eg{{t}_{1}}+{{g}^{2}}{{t}_{0}}$ and ${{t}_{6}}=e{{t}_{1}}-2g{{t}_{0}}$, into the above formula, we obtain
\begin{eqnarray} \label{F1}
\mathbf{F}_1(\xi)=\underbrace{\log ({{\Xi}_{1}})\left[ \frac{{{\Xi}_{0}}\log ({{\Xi}_{0}}){{t}_{6}}-{{t}_{3}}}{{{e}^{3}}{{\Xi}_{0}}}+\frac{\xi{{t}_{0}}}{{{e}^{2}}} \right]}_{\mathbf{A}}
+\underbrace{\frac{t_6}{{{e}^{3}}}\int{\frac{({{t}_{1}}+2\xi{{t}_{0}})\log (\Xi_0)}{\Xi_1}d\xi}}_{\mathbf{B}_1}+\underbrace{\int{\frac{({{t}_{1}}+2\xi{{t}_{0}})\left[ -t_3 + e \xi \Xi_0 t_0 \right]}{e^3 \Xi_0 \Xi_1}d\xi}}_{\mathbf{B}_2},
\end{eqnarray}
where if $t_{1}^{2}-4{{t}_{0}}{{t}_{2}}\ge 0$,
\begin{equation}\label{B1_0}
\mathbf{B}_1 = \frac{t_6}{e^3} \bigintsss \log (\Xi_0)\left( \frac{2{{t}_{0}}}{2{{t}_{0}}\, \xi+\sqrt{t_{1}^{2}-4{{t}_{0}}{{t}_{2}}}+{{t}_{1}}}+\frac{2{{t}_{0}}}{2{{t}_{0}}\, \xi-\sqrt{t_{1}^{2}-4{{t}_{0}}{{t}_{2}}}-{{t}_{1}}} \right) d\xi,
\end{equation}
and otherwise
\begin{equation}\label{B1_1}
\mathbf{B}_1 = \frac{t_6}{e^3} \bigintsss \log (\Xi_0)\left( \frac{2{{t}_{0}}}{2{{t}_{0}}\, \xi+i\sqrt{-t_{1}^{2}+4{{t}_{0}}{{t}_{2}}}+{{t}_{1}}}+\frac{2{{t}_{0}}}{2{{t}_{0}}\, \xi-i\sqrt{-t_{1}^{2}+4{{t}_{0}}-{{t}_{1}}}{2{{t}_{0}}}} \right) d\xi.
\end{equation}
The integration of Eq. (\ref{B1_0}) results in
\begin{eqnarray}\label{B1_2}
\mathbf{B}_1&=&-\frac{t_6}{{{e}^{3}}} \left[\log ({{\Xi}_{2}}+s)\log \left( -\frac{2{{\Xi}_{0}}{{t}_{0}}}{s_1} \right) + \log \left( {{\Xi}_{2}}-s \right)\log \left( -\frac{2{{\Xi}_{0}}{{t}_{0}}}{s_2} \right)  + \text{Li}_2\left( \frac{e\Xi_2+es}{s_1} \right) +  \text{Li}_2\left( \frac{e{{\Xi}_{2}}-es}{s_2} \right) \right.  \\ \nonumber
 &-& \left. \log ({{\Xi}_{0}})\log (4{{\Xi}_{1}}{{t}_{0}}) \right],
\end{eqnarray}
with $\Xi_2=t_1+2\xi t_0$, $s=\sqrt{\left| t_{1}^{2}-4{{t}_{0}}{{t}_{2}} \right|}$, $s_1=t_6+es$ and $s_2=t_6-es$. The dilogarithm (Spence's function) is defined as $\text{Li}_2(z)=-\int_0^z{\log(1-t)dt}$. Eq. (\ref{B1_1}) becomes
\begin{eqnarray}
\mathbf{B}_1&=&\frac{ t_6}{{{e}^{3}}}\left\{\log ({{\Xi}_{0}})\left[ \log (\frac{{{\Xi}_{1}}}{{{t}_{0}}})+2\log (2{{t}_{0}})\right] -2\operatorname{atan2}(s,\Xi_2)\operatorname{atan2}(es,t_6)  -2 \text{Li}_2\left[ \frac{e(2e{{t}_{2}}+e\xi{{t}_{1}}-g\Xi_2)+ie{{\Xi}_{0}}s}{2{{t}_{0}}{{t}_{3}}} \right] \right.
 \nonumber \\
&-&\left. \log (4{{t}_{0}}{{\Xi}_{1}})\log \left( \frac{ {{\Xi}_{0}} \sqrt{{{t}_{0}}{{t}_{3}}}}{ {{t}_{3}} } \right) \right\}.
\end{eqnarray}
Similarly, in Eq. (\ref{F1}) we have $\mathbf{B}_2 = \frac{1}{{{e}^{3}}} \left[  -\log ({{\Xi}_{0}}){{t}_{6}}-g\log ({{\Xi}_{1}}){{t}_{0}}+2e\xi{{t}_{0}}+es\, \log \left( \frac{{{\Xi}_{2}}-s}{{{\Xi}_{2}}+s} \right) \right]$ if $t_{1}^{2}-4{{t}_{0}}{{t}_{2}}\ge 0$, otherwise $\mathbf{B}_2 = \frac{1}{{{e}^{3}}} \left[  -\log ({{\Xi}_{0}}){{t}_{6}}-g\log ({{\Xi}_{1}}){{t}_{0}}+2e\xi{{t}_{0}}-2es\, \text{atan}(\frac{\Xi_2}{s}) \right]$. 

\begin{figure}
\begin{center}
\includegraphics[width=0.7\textwidth]{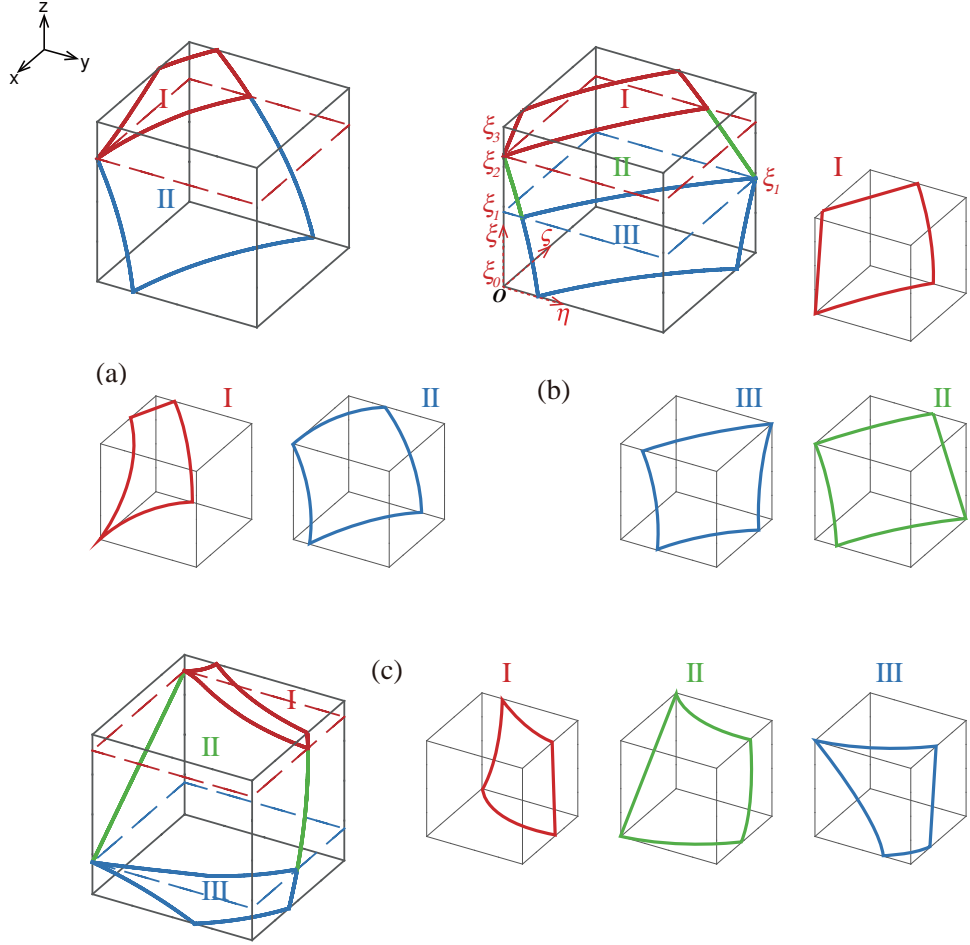}
\caption{The strategy to split the 3D complicated cases into elementary cases in Fig. \ref{3d_basic}.}
\label{3d_multi}
\end{center}
\end{figure}

To calculate $\mathbf{F}_2(\xi)$, we define $\mathbf{F}_2(\xi)=\mathbf{C}+\mathbf{D}+\mathbf{E}$, with $\mathbf{C}$, $\mathbf{D}$ and $\mathbf{E}$ being
\begin{eqnarray}\label{C}
\mathbf{C}&=&\int -\frac{\log \left( f\xi+h \right)\left( {{t}_{2}}+\xi{{t}_{1}}+{{\xi}^{2}}{{t}_{0}} \right)}{{{(e\xi+g)}^{2}}} d\xi = -\log ({{\Xi}_{3}})\frac{\log ({{\Xi}_{0}}){{t}_{6}}}{{{e}^{3}}} -\log ( {{\Xi}_{3}})\left( -\frac{{{t}_{3}}}{{{e}^{3}}{{\Xi}_{0}}}+\frac{\xi{{t}_{0}}}{{{e}^{2}}} \right)-\frac{f\log ( {{\Xi}_{3}}){{t}_{3}}}{{{e}^{3}}{{t}_{4}}}
\nonumber \\
&+&\frac{\xi{{t}_{0}}}{{{e}^{2}}} + \frac{{{t}_{6}}}{{{e}^{3}}} \log ({{\Xi}_{0}})\log (\frac{e{{\Xi}_{3}}}{{{t}_{4}}})+\frac{{{t}_{6}}}{{{e}^{3}}}\text{Li}_2\left( \frac{-f{{\Xi}_{0}}}{{{t}_{4}}} \right) + \frac{\log ( {{\Xi}_{3}})\left( {{f}^{2}}{{t}_{3}}-eh{{t}_{0}}{{t}_{4}} \right)}{{{e}^{3}}f{{t}_{4}}},
\end{eqnarray}
\begin{eqnarray}\label{D}
\mathbf{D}&=&\int -\frac{\log \left( a\xi+c \right)\left( {{t}_{2}}+\xi{{t}_{1}}+{{\xi}^{2}}{{t}_{0}} \right)}{{{(e\xi+g)}^{2}}} d\xi = -\frac{{{t}_{6}}}{{{e}^{3}}}\log ({{\Xi}_{4}})\log ({{\Xi}_{0}})-\log ( {{\Xi}_{4}})\left( -\frac{{{t}_{3}}}{{{e}^{3}}{{\Xi}_{0}}}+\frac{\xi{{t}_{0}}}{{{e}^{2}}} \right) -\frac{a\log ( {{\Xi}_{0}} ){{t}_{3}}}{{{e}^{3}}{{t}_{5}}}
\nonumber \\
&+&\frac{\xi{{t}_{0}}}{{{e}^{2}}} + \frac{{{t}_{6}}}{{{e}^{3}}}\log ({{\Xi}_{0}})\log (\frac{e{{\Xi}_{4}}}{{{t}_{5}}})+ \frac{{{t}_{6}}}{{{e}^{3}}}\text{Li}_2\left( \frac{- a{{\Xi}_{0}}}{{{t}_{5}}} \right)
+\frac{\log ({{\Xi}_{4}})\left( {{ a}^{2}}{{t}_{3}}-ec{{t}_{0}}{{t}_{5}} \right)}{{{e}^{3}}a {{t}_{5}}}
\end{eqnarray}
and
\begin{eqnarray}\label{E}
\mathbf{E}=\int{\frac{b\xi+d}{e\xi+g}d\xi}=\frac{b\xi}{e}-\frac{(bg-de)\log ( \Xi_0)}{{{e}^{2}}},
\end{eqnarray}
respectively, with ${{t}_{4}}=eh-fg$, ${{t}_{5}}=ce-ag$, ${{\Xi}_{3}}=f\xi+h$ and ${{\Xi}_{4}}=a\xi+c$. We can generalize Eq. (\ref{volume3d}) for cases with different $i_{\zeta}$,
\begin{eqnarray}\label{eqalpha}
\alpha &=&  (1-2i_{\varsigma})\left[ \mathbf{F}\left( \xi_{1}  \right)^{*}- \mathbf{F}\left( \xi_{0}  \right)^{*}\right]
+i_{\varsigma}\left[ \mathbf{G}\left( \xi_{1}  \right)^{*}- \mathbf{G}\left( \xi_{0}  \right)^{*}\right],
\end{eqnarray}
where $\mathbf{F}^{*}\left( \xi  \right)=\left( 1-2{{i}_{\eta }} \right)\left( \mathbf{A^{*}}-{{\mathbf{B}}_{\mathbf{1}}}^{*}-{{\mathbf{B}}_{\mathbf{2}}}^{*} \right)+\left( 1-2{{i}_{\eta }} \right)\left( \mathbf{C}^{*}+\mathbf{D}^{*}+\mathbf{E}^{*} \right)+\mathbf{G}^{*}$ and $\mathbf{G}\left( \xi  \right)^{*} =\left( 1-2{{i}_{\eta }} \right)\mathbf{G}\left( \xi  \right)$, where ${{\mathbf{A}}}^{*}$, ${{\mathbf{B}}_{\mathbf{1}}}^{*}$ and ${{\mathbf{B}}_{\mathbf{2}}}^{*}$ are the same as ${{\mathbf{A}}}$, ${{\mathbf{B}}_{\mathbf{1}}}$ and ${{\mathbf{B}}_{\mathbf{2}}}$, respectively. ${\mathbf{C}^{*}}$, ${\mathbf{D}^{*}}$, ${\mathbf{E}^{*}}$ and ${\mathbf{G}^{*}}$ are defined as
\begin{eqnarray}
{\mathbf{C}^{*}}&=&-\frac{{{t}_{6}}}{{{e}^{3}}} \log ({{\Xi}_{0}}) \left[ \log ({{\Xi}_{3}}+{{i}_{\eta }}{{\Xi}_{0}})  - \log \frac{e\left( {{\Xi}_{3}}+{{i}_{\eta }}{{\Xi}_{0}} \right)}{{{t}_{4}}} \right] + \log ( {{\Xi}_{3}}+{{i}_{\eta }}{{\Xi}_{0}} ) \left[ \frac{ {{\left( f+{{i}_{\eta }}e \right)}^{2}}{{t}_{3}}-e\left( h+{{i}_{\eta }}e \right){{t}_{0}}{{t}_{4}} }{{{e}^{3}}\left( f+{{i}_{\eta }}e \right){{t}_{4}}} + \frac{{{t}_{3}}}{{{e}^{3}}{{\Xi}_{0}}}-\frac{\xi{{t}_{0}}}{{{e}^{2}}}\right]
\nonumber \\
&+&\frac{{{t}_{6}}}{{{e}^{3}}}\text{Li}_2\left[ \frac{-\left( f+{{i}_{\eta }} \right)e{{\Xi}_{0}}}{{{t}_{4}}} \right] -\frac{\left( f+{{i}_{\eta }}e \right)\log ( {{\Xi}_{0}}){{t}_{3}}}{{{e}^{3}}{{t}_{4}}}+\frac{\xi{{t}_{0}}}{{{e}^{2}}},
\end{eqnarray}
\begin{eqnarray}
{\mathbf{D}^{*}}&=&-\frac{{{t}_{6}}}{{{e}^{3}}}\log ({{\Xi}_{0}}) \left[ \log ({{\Xi}_{4}}+{{i}_{\varsigma }}{{\Xi}_{0}})- \log \frac{e\left( {{\Xi}_{4}}+{{i}_{\varsigma }}{{\Xi}_{0}} \right)}{{{t}_{5}}}\right] + \log ({{\Xi}_{4}}+{{i}_{\varsigma }}{{\Xi}_{0}}) \left[ \frac{{{\left( a+{{i}_{\varsigma }}e \right)}^{2}}{{t}_{3}}-e\left( c+{{i}_{\varsigma }}g \right){{t}_{0}}{{t}_{5}} }{{{e}^{3}}\left( a+{{i}_{\varsigma }}e \right){{t}_{5}}}+\frac{{{t}_{3}}}{{{e}^{3}}{{\Xi}_{0}}}-\frac{\xi{{t}_{0}}}{{{e}^{2}}}\right]
\nonumber \\
&+&\frac{{{t}_{6}}}{{{e}^{3}}}\text{Li}_2\left[ \frac{-\left( a+{{i}_{\varsigma }}e \right){{\Xi}_{0}}}{{{t}_{5}}} \right]-\frac{\left( a+{{i}_{\varsigma }}e \right)\log ( {{\Xi}_{0}}){{t}_{3}}}{{{e}^{3}}{{t}_{5}}}+\frac{\xi{{t}_{0}}}{{{e}^{2}}},
\end{eqnarray}
\begin{eqnarray}
{\mathbf{E}^{*}}&=&-{{i}_{\varsigma }}\frac{(eh+ah-fg-bg-cf+de+ad-bc)}{{{(a+e)}^{2}}}\log ( {{\Xi}_{0}}+{{\Xi}_{4}})
-\frac{(bg-de)+{{i}_{\eta }}(ag-ce)-{{i}_{\varsigma }}(eh-fg)}{{{e}^{2}}}\log ({{\Xi}_{0}})
\nonumber \\
&+&\frac{\left( 1-{{i}_{\varsigma }}+{{i}_{\varsigma }}a \right)({{i}_{\eta }}{{i}_{\varsigma }}e+{{i}_{\varsigma }}f+{{i}_{\eta }}a+b)\xi}{e\left[ 1-{{i}_{\varsigma }}+(a+e){{i}_{\varsigma }} \right]},
\end{eqnarray}
and ${\mathbf{G}^{*}}=\left( 1-2{{i}_{\eta }} \right){\mathbf{G}}+\xi{{i}_{\eta }}$, respectively. Now, we can use this general formulation to calculate 3D volume fraction of all cases sketched in Fig. \ref{3d_basic}. These elementary cases generate more complex cases as indicated in Fig. \ref{3d_multi}, and are defined as type I, II, and III, respectively. These complex cases can be split them into elementary cases in Fig. \ref{3d_basic}, as illustrated in each subgraph. For instance, the case in Fig. \ref{3d_multi}(a) is composed of type I and II, whereas cases in Fig. \ref{3d_multi}(b) and (c) are combinations of all three elementary types. In order to compute the volumes of those cases, the first step is to decompose the cell into elementary types, which is easy to accomplish. For example, in Fig. \ref{3d_multi}(b), along the $\xi$-direction we can find two intersection points $\xi_{1}$ and $\xi_{2}$ which occur at two different edges parallel to $\xi$-direction. By these intersection points, the whole cell is cut into three parts in $\xi$ directions, as shown in Fig. \ref{3d_multi}(b).

\subsection{Accuracy test}
\begin{figure}
\begin{center}
\subfloat[][]{\label{accuracy2da}\includegraphics[scale=0.25]{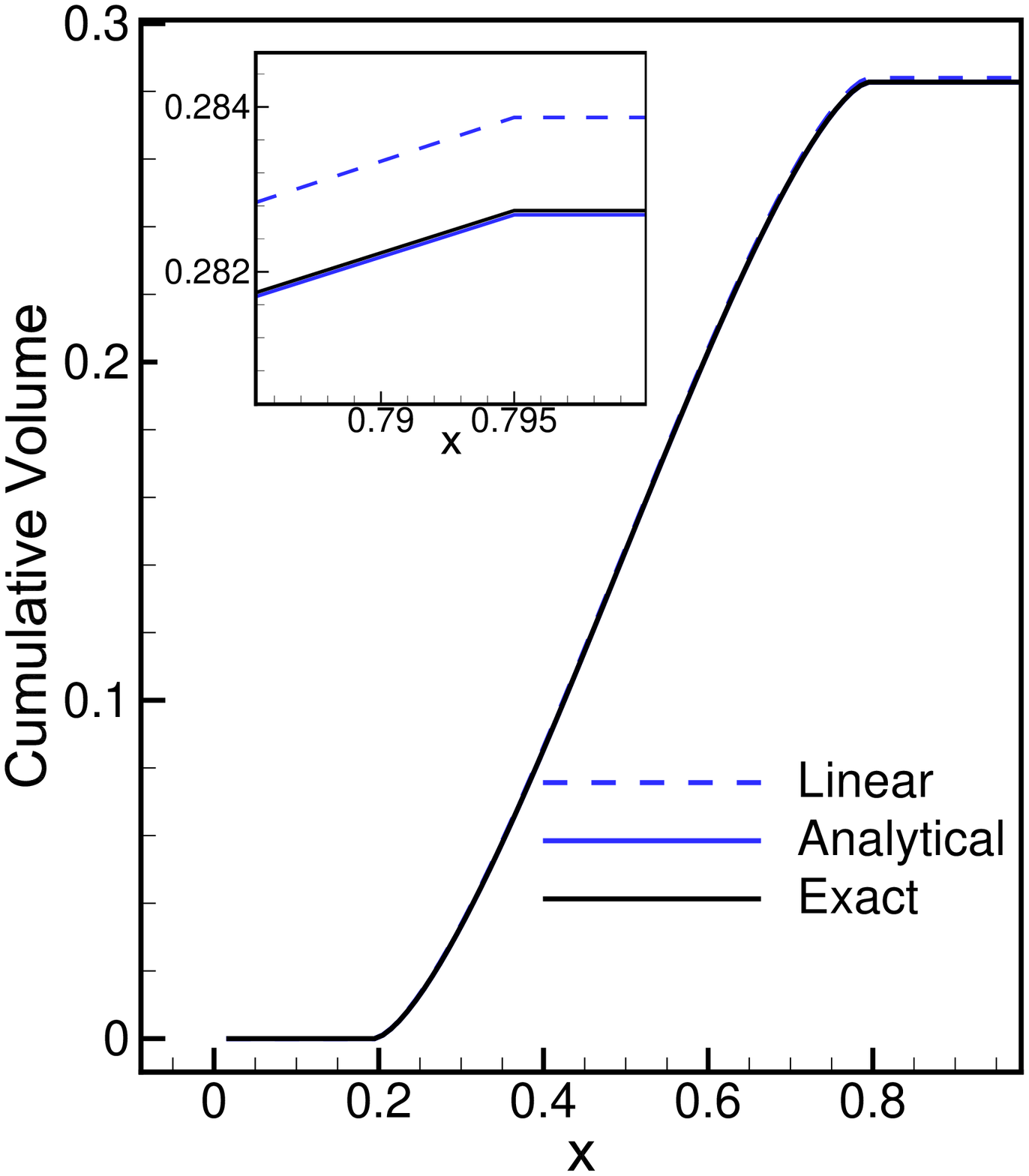}}
\subfloat[][]{\label{accuracy2db}\includegraphics[scale=0.25]{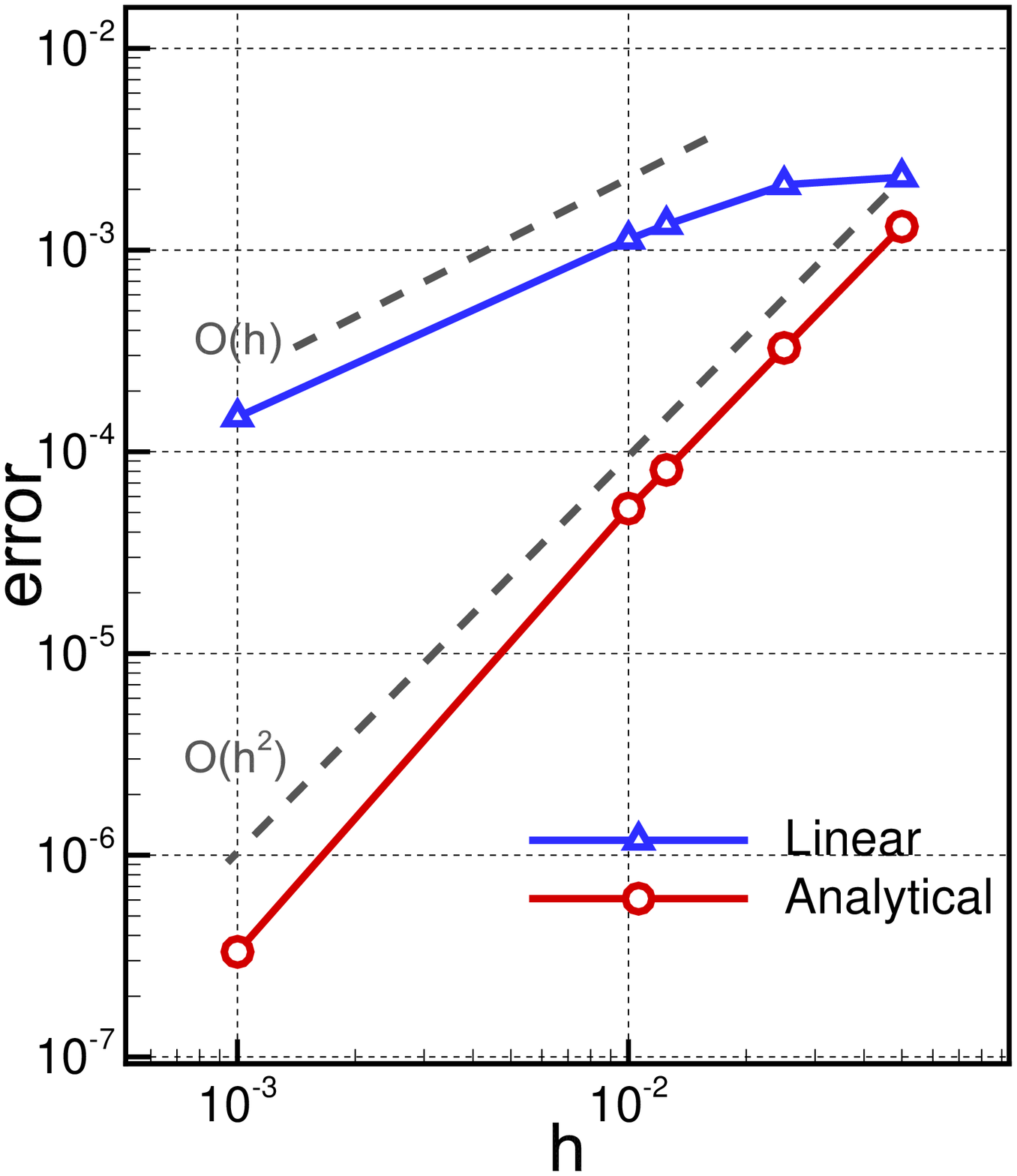}}
\subfloat[][]{\label{accuracy2dc}\includegraphics[scale=0.25]{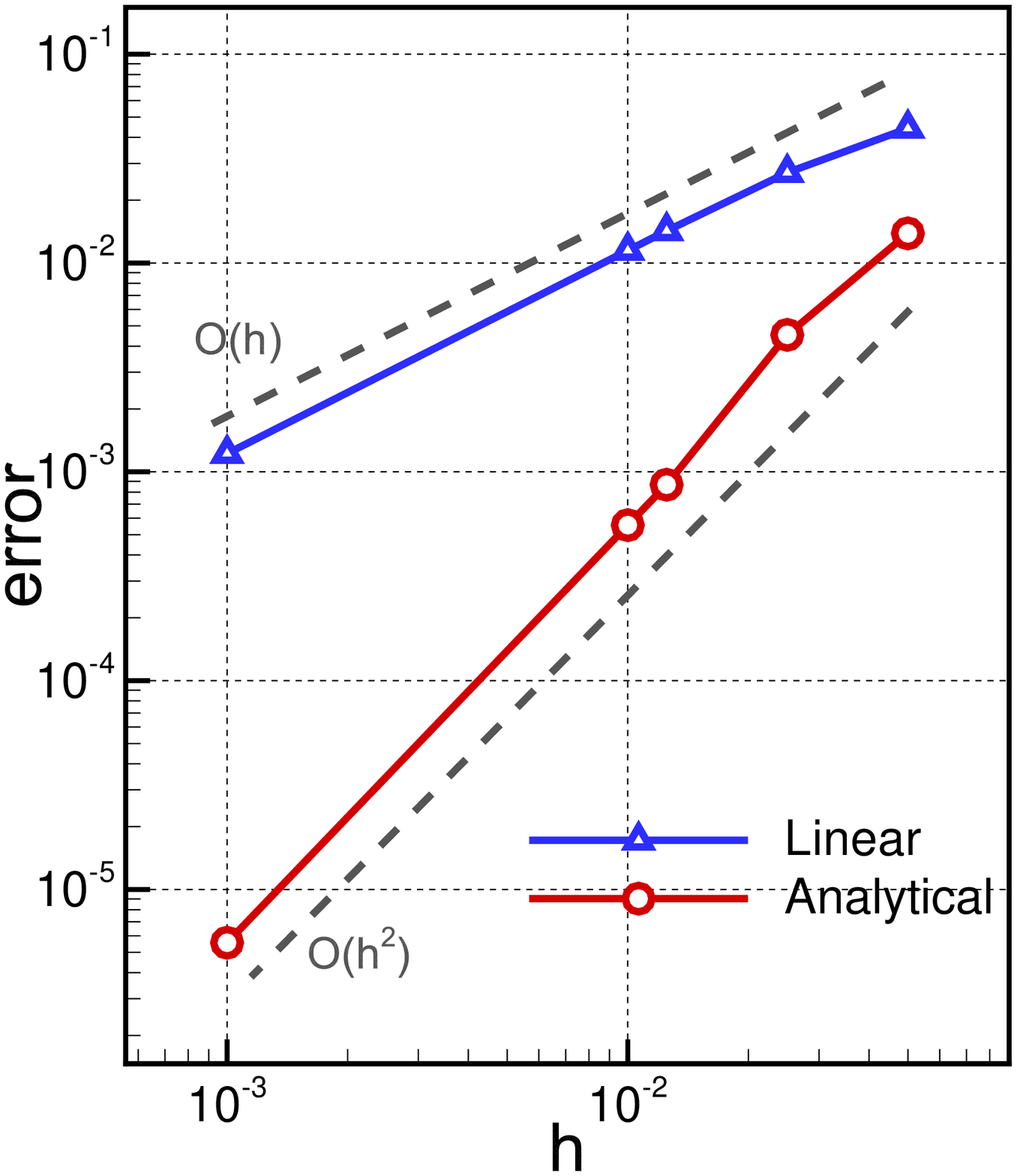}}
\caption{\label{accuracy2d} The accuracy test for 2D and 3D simple cases: (a) comparison of 2D cumulative volumes, (b) convergence for the 2D absolute errors and (c) convergence for the 3D absolute errors.}
\end{center}
\end{figure}

We calculate the volumes of a circle and a sphere to assess the accuracy of our analytic formulations in 2D and 3D. The results of our method are compared to those of a linear approximation method \cite{lauer2012numerical}, as plotted in Fig. \ref{accuracy2d}.

For 2D, as shown in Fig. \ref{accuracy2d}(a), the cumulative volumes along $x$-direction indicate that both linear and analytical results are in good agreement with the exact result while locally the analytical formulation is more accurate. In Fig. \ref{accuracy2d}(b) 2nd-order convergence rate is demonstrated for our formulation because it is based on bilinear/triliner interpolation while the linear method is 1st-order as expected. The magnitude of the error is much smaller in our method. The result of the 3D case in Fig. \ref{accuracy2d}(c) exhibits similar error distribution with Fig. \ref{accuracy2db}. As a consequence, for 2D and 3D, our analytical formulation provides more accurate volume-estimation than the linear method due to a more accurate interface representation.

\subsection{Consistency test}
\vspace{0.2cm}

\noindent \textbf{Definition.} The consistency condition of the volume-estimation during mesh refinement is $\int_{\Omega} \alpha^{0} d V  = \int_{\Omega} \alpha^{\ell} d V$, i.e., volumes calculated from the initial volumetric field and the refined field are exactly the same.

\vspace{0.2cm}

The significance of maintaining consistency of volume estimation is obvious. For example, for an adaptive multiphase solver \cite{sussman1999adaptive} conservation errors can only be prevented if the volumes have an identical value at each resolution level. This can not be achieved by a linear representation. All such numerical volume estimations are not consistent as the errors introduced by numerical methods depend on grid resolution. 

\vspace{0.4cm}
\noindent \textbf{Theorem.} The analytical volume-estimation formulation is consistent if the level-set field is refined by bilinear/trilinear interpolation.
\begin{pf}
The necessary and sufficient condition for consistency is that the parameters of Eq. (\ref{eqvar2dphi}) or Eq. (\ref{eq3d}) remain invariant upon mesh refinement as no numerical error is introduced by the analytical formulation and the coordinate transformation doest not change the computed volume. This can be easily proved. For simplicity, we consider the 2D cases for description. Assume the grid point level-set values $\phi^{\ell}$ of a cell $\left[ i,i+1 \right] \times \left[ j,j+1 \right] $ on level $\ell$ are given. During refinement, a cell is divided into $4$ subcells. The level-set values $\phi^{\ell+1}$ at level $\ell+1$ are computed by bilinear interpolation. Thus the parameters of the interface equation for the subcell $\left[ 2i+1,2i+2 \right] \times \left[ 2j+1,2j+2 \right]$ are 
\begin{equation}
\beta^{'}_0=\phi^{\ell+1}_{00}, \quad  \beta^{'}_1=\phi^{\ell+1}_{10}-\phi^{\ell+1}_{00},\quad  \beta^{'}_2=\phi^{\ell+1}_{01}-\phi^{\ell+1}_{00}, \quad \beta^{'}_3=\phi^{\ell+1}_{00}+\phi^{\ell+1}_{11}-\phi^{\ell+1}_{01}-\phi^{\ell+1}_{10}
\end{equation}
which can be expressed by $\phi^{\ell}$ with bilinear interpolation,
\begin{eqnarray}
&&\beta^{'}_0=\frac{\phi^{\ell}_{00}+\phi^{\ell}_{01}+\phi^{\ell}_{10}+\phi^{\ell}_{11}}{4}, \,  \beta^{'}_1=\frac{\phi^{\ell}_{10}+\phi^{\ell}_{11}}{2}-\frac{\phi^{\ell}_{00}+\phi^{\ell}_{01}+\phi^{\ell}_{10}+\phi^{\ell}_{11}}{4},\, \beta^{'}_2 = \frac{\phi^{\ell}_{01}+\phi^{\ell}_{11}}{2}-\frac{\phi^{\ell}_{00}+\phi^{\ell}_{01}+\phi^{\ell}_{10}+\phi^{\ell}_{11}}{4} \\
\nonumber
&&\beta^{'}_3 =\frac{\phi^{\ell}_{00}+\phi^{\ell}_{01}+\phi^{\ell}_{10}+\phi^{\ell}_{11}}{4} + \phi^{\ell}_{11}-\frac{\phi^{\ell}_{01}+\phi^{\ell}_{11}}{2}-\frac{\phi^{\ell}_{10}+\phi^{\ell}_{11}}{2}.
\end{eqnarray}
This implies that the interface equation at $\ell+1$,
\begin{equation}
\phi(x',y')=\beta^{'}_0+\beta^{'}_1 x'+\beta^{'}_2 y'+\beta^{'}_3 x'y'=0,
\end{equation}
is the same as Eq. (\ref{eqvar2dphi}), as $x' = 2\, x-1$ and $y' = 2\, y-1$. For the other subcell we can obtain the same conclusion, which completes the proof.  $\square$
\end{pf}
For example, a single cell test with specific level-set values clearly shows the consistency feature of our method.
After prescribing the grid values, $\phi^{\ell}_{00}$ = 0.1, $\phi^{\ell}_{10}$ = 0.6, $\phi^{\ell}_{01}$ = -0.3, and $\phi^{\ell}_{11}$ = -0.1, we obtain the coefficients $a=0.5$, $b=0.1$, $c=-0.3$, and $d=-0.4$. The origin is placed at $(0,0)$ and the integration variable $\eta$ is $x$, hence the interface profile becomes $\zeta=(5\eta+1)/(3\eta+4)$ with an integration range of $[0, 1]$. After substituting the parameters into the Eq. (\ref{volume2d}), the volume of this piece on this level is $\alpha=\frac{1}{9}(17\, \mathrm{log}(4)-17\, \mathrm{log}(7) + 15)$. The bilinear interpolation leads to $\phi^{\ell+1}_{00}$ = 0.1, $\phi^{\ell+1}_{10}$ = 0.35, $\phi^{\ell+1}_{01}$ = -0.1, and $\phi^{\ell+1}_{11}$ = 0.075 for subcell cell $\left[ 2i,2i+1 \right] \times \left[ 2j,2j+1 \right]$. Substitute these values into Eq. (\ref{volume2d}), we can get a volume of $\alpha_{00}=\frac{1}{63} (120-476\,\mathrm{log}\left( \frac{68}{7}\right)) +\frac{68}{9}\,\mathrm{log}\left( 8\right) +\frac{3}{7} $. Analogously, the volumes of the other three subcells are $\alpha_{10}=1.0$, $\alpha_{01}=\frac{68}{9}(\,\mathrm{log}\left( \frac{68}{7}\right) -12 - 68\,\mathrm{log}\left( 11\right) + 21)$ and $\alpha_{11}=\frac{68}{9}(\,\mathrm{log}\left( 11\right) - \,\mathrm{log}\left( 14\right) + 21)$, respectively. Obviously, the consistency is satisfied as $\alpha=(\alpha_{11}+\alpha_{10}+\alpha_{01}+\alpha_{00})/4$. Analogously, we can proof the consistency of our volume-estimation method during coarsening.

\begin{figure}
\begin{center}
\subfloat[][]{\label{consistancya} \includegraphics[scale=0.25]{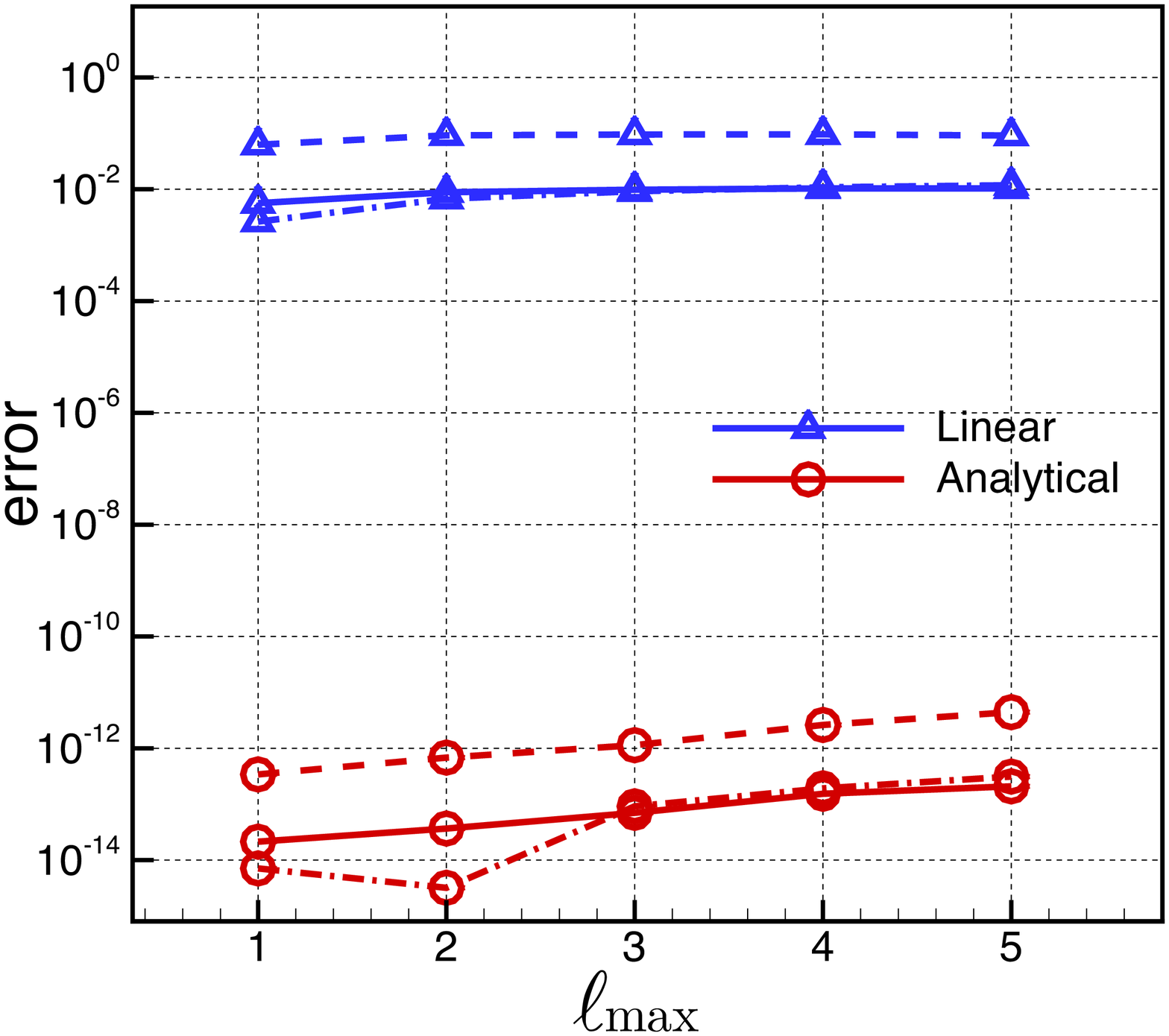}}
\subfloat[][]{\label{consistancyb} \includegraphics[scale=0.25]{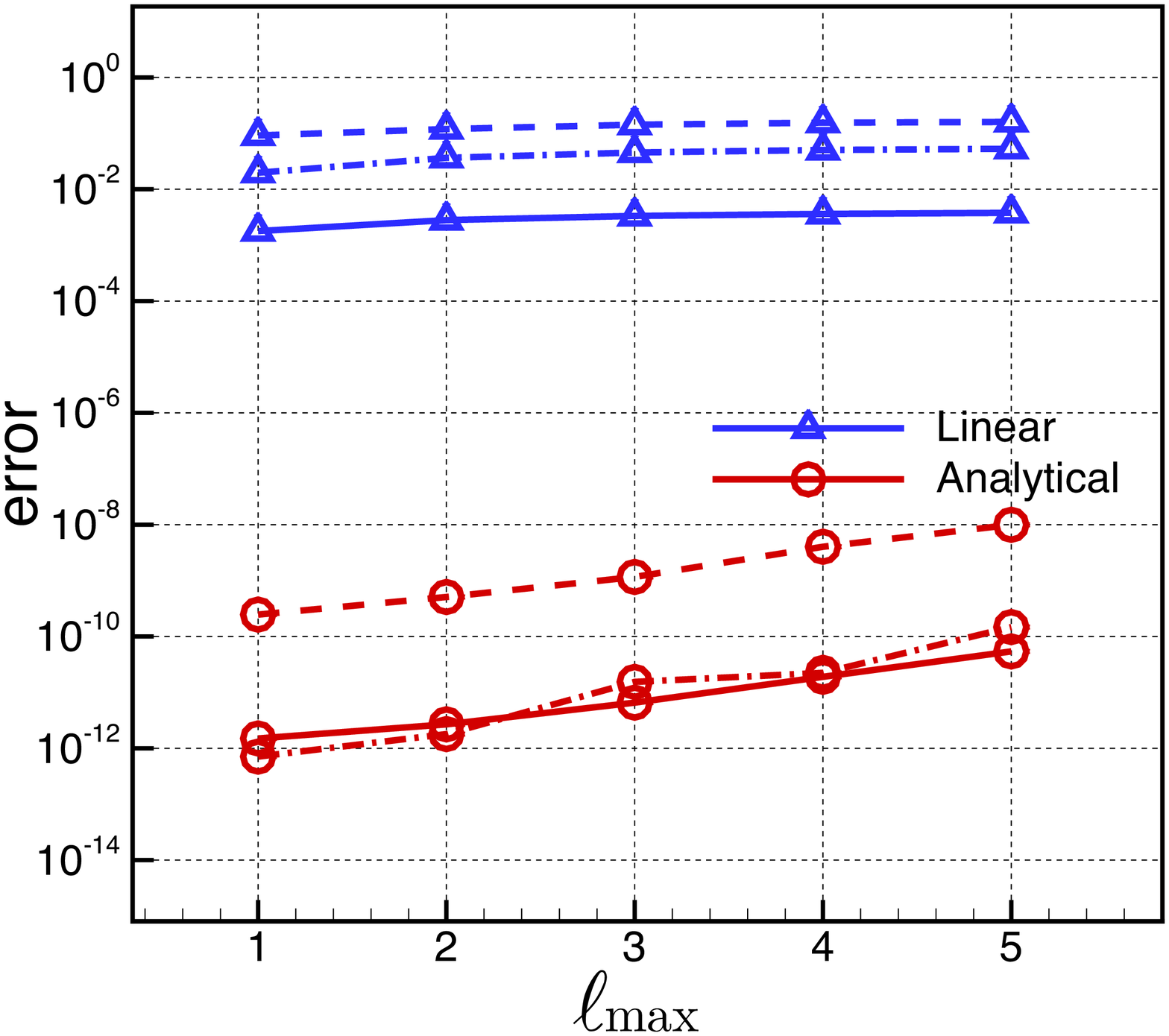}}
\caption{\label{consistancy} 2D consistency for $L_1$(dash-dotted line), $L_2$(solid line) and $L_\infty$(dashed line). The accuracy test: (a) 2D and (b) 3D.}
\end{center}
\end{figure}

In addition, simple cases used in the accuracy test are employed to test the consistency of the present method. As shown in Fig. \ref{consistancy}, the initial resolution is $h=0.05$ on $\ell = 0$, then the level-set fields are refined from $\ell_{max}=1$ to $5$. The dash-dotted, solid and dashed lines stand for error norms,
\begin{equation}\label{eq_error}
L_1 = \frac{\sum_{i,j}|\alpha^{\ell}_{i,j}-\alpha^{0}_{i,j}|}{\sum_{i,j} \alpha^{\ell}_{i,j}}, \quad L_2=[\frac{\sum_{i,j}|\alpha^{\ell}_{i,j}-\alpha^{0}_{i,j}|^2}{\sum_{i,j} \alpha^{\ell}_{i,j}}]^{\frac{1}{2}} \quad \mathrm{and} \quad L_\infty=\max_{i,j}|\alpha^{\ell}_{i,j}-\alpha^{0}_{i,j}|,
\end{equation}
respectively, with the superscripts `$0$' and `$\ell$' being the volume fraction of the cell $[i,j]\times[i+1,j+1]$ on $\ell=0$ and the summation of volume fractions on $\ell>0$ within the same cell. 
In 2D, with $\ell_{\text{max}}$ increasing, the discrepancy between $\ell = 0$ and $\ell_{\text{max}}$ for linear approximation method is large, ranging from $5\times10^{-2}$ to $1\times10^{-1}$ for $L_\infty$. The linear reconstruction inside each cell is responsible for this distinct errors as it does not have a consistent representation of interface between different resolutions.
With respect to the analytical formulation the magnitudes are much lower than those of the linear method. So the bilinear/trilinear representation and interpolation are beneficial to preserving the consistency of the volume estimation. The 2D and 3D errors produced in our method can be neglected. Therefore it is concluded that for both 2D and 3D our method automatically achieves consistency which is attributed to the analytical formulation and the bilinear/trilinear interpolation. 

\subsection{Additional remarks}
\begin{enumerate}
\item Note that only resolved cases have an analytical formulation. If there exists ambiguity such as multiple intersections in one cell edge, pure analytical method may not provide the correct results as the formulation is not able to resolve this cell. One can easily address this issue by adaptive refinement in the ambiguity regions or performing an interface scale separation \cite{han2015scale}.
\item Our method achieves 2nd-order accuracy which is generally acceptable for most of the applications. However the goal of our work is not to achieve high-order accuracy, but to maintain the consistency which is the unique feature of our method.
\end{enumerate}

\section{Applications} \label{sec:app}

\begin{figure}
\begin{center}
\subfloat[][]{\label{geo_resulta} \includegraphics[scale=0.25]{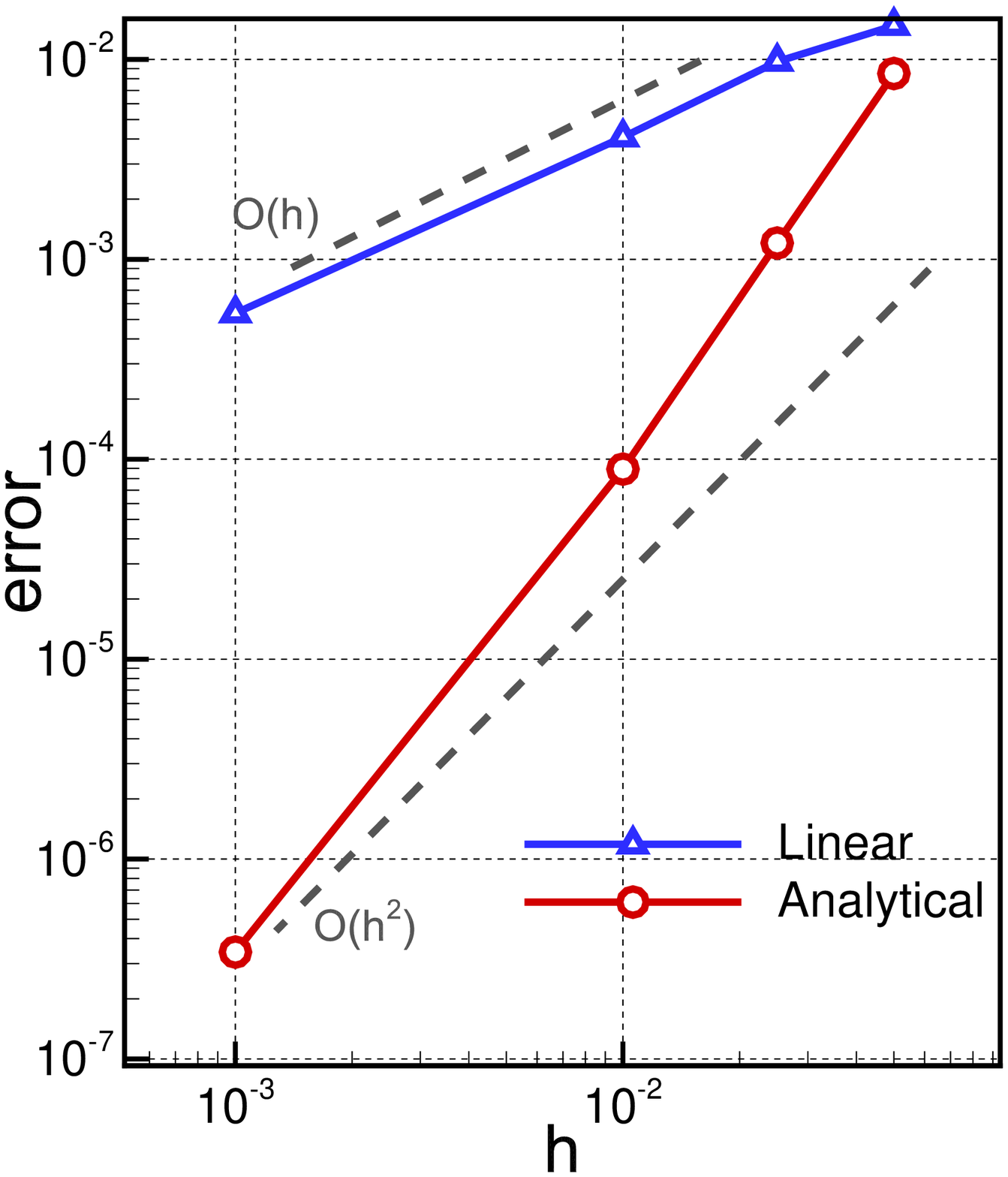}}
\subfloat[][]{\label{geo_resultb} \includegraphics[scale=0.25]{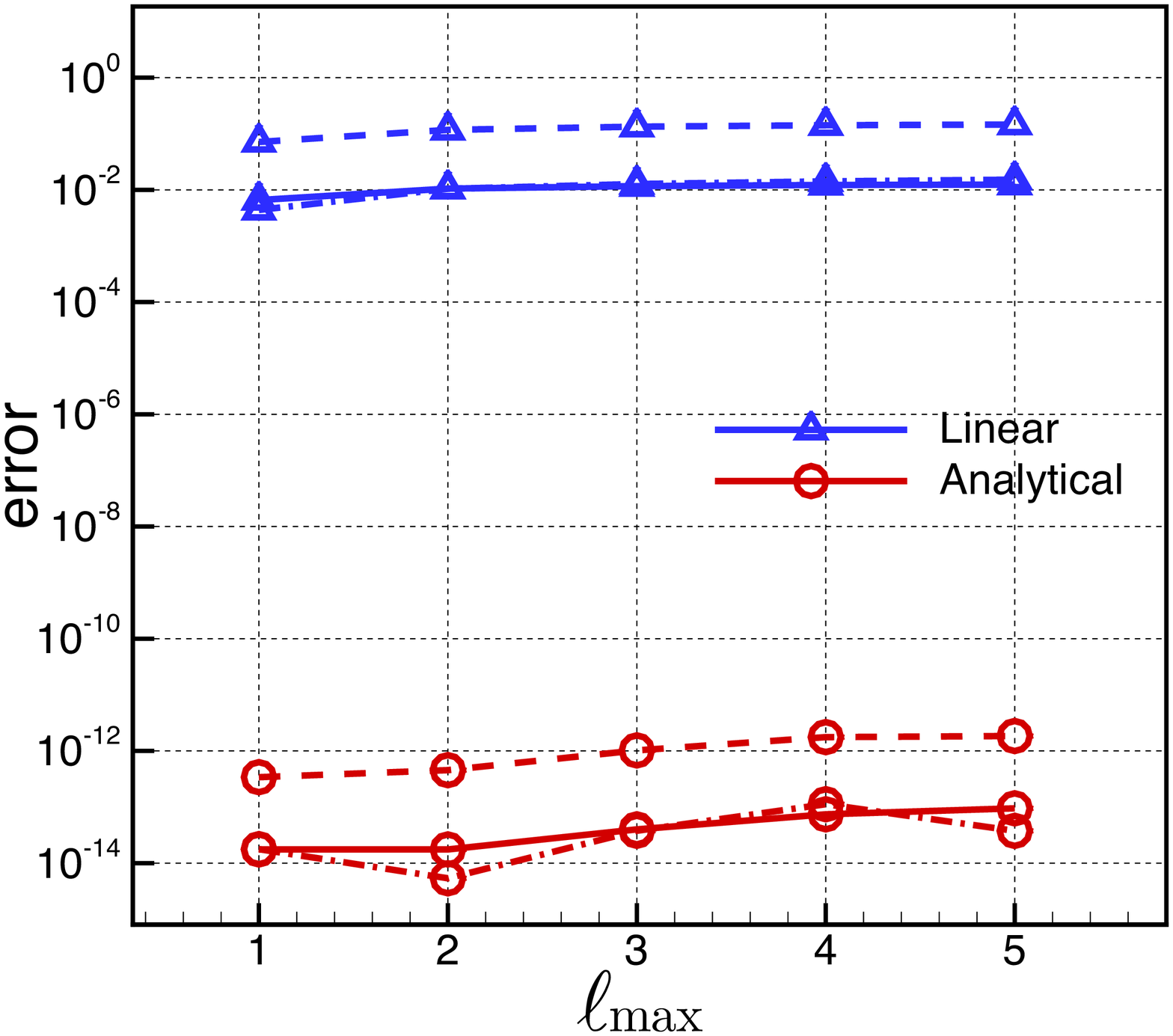}} \\
\subfloat[][]{\label{geo_resultc} \includegraphics[scale=0.25]{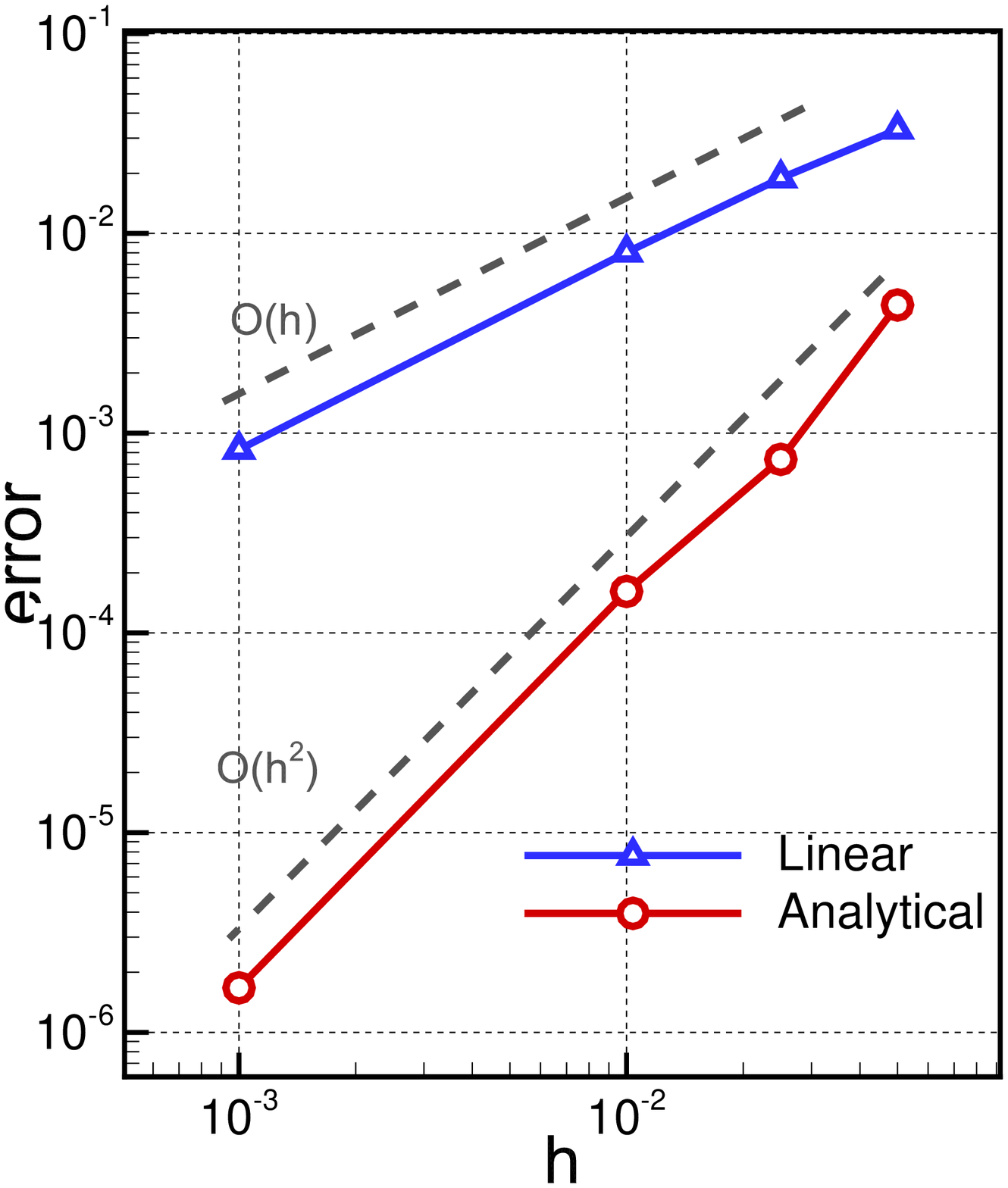}}
\subfloat[][]{\label{geo_resultd} \includegraphics[scale=0.25]{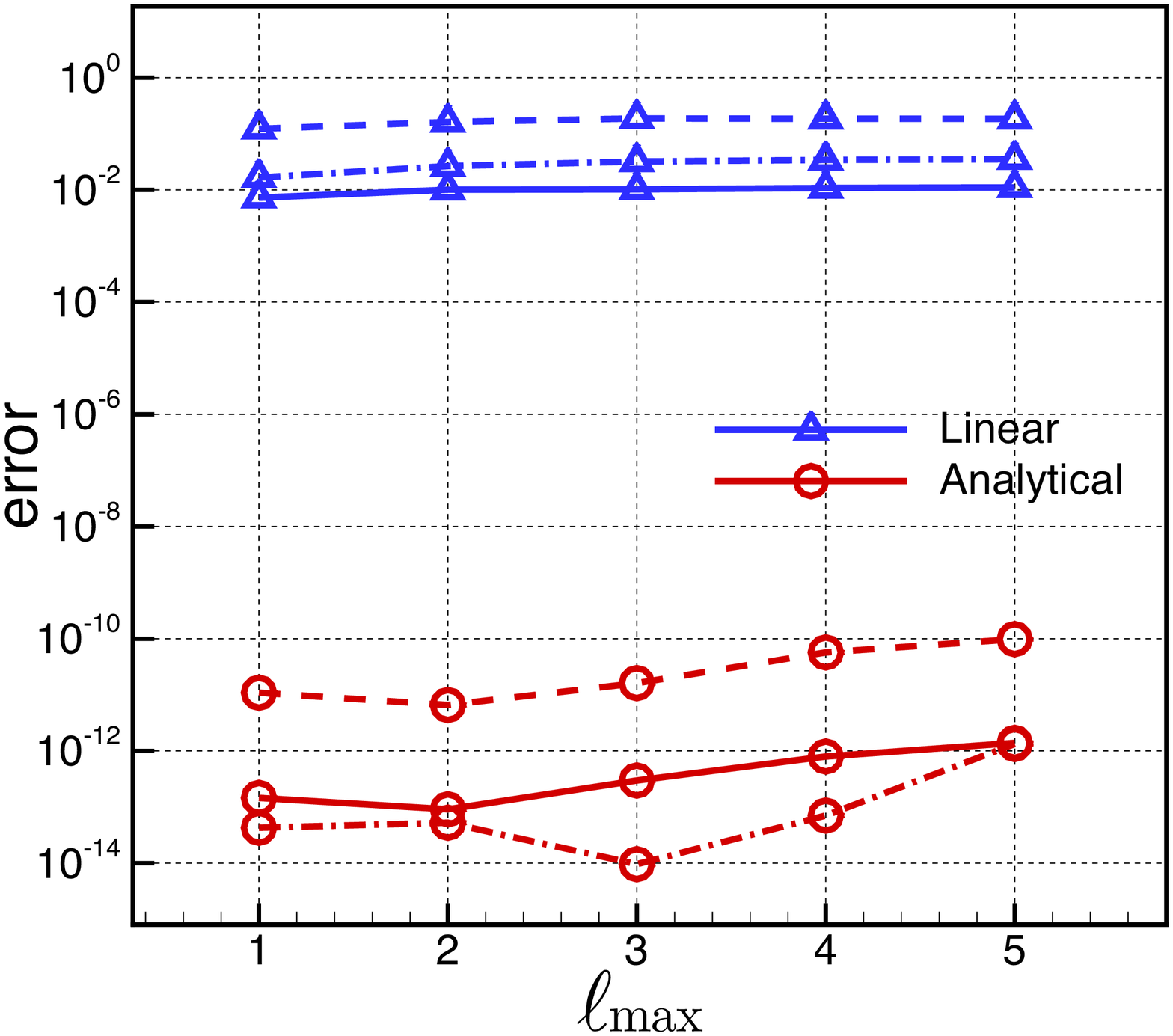}}
\caption{\label{geo_result} Accuracy and consistency tests for double circles (a and b) and Zalesak disk (c and d).}
\end{center}
\end{figure}

\begin{figure}
\begin{center}
\subfloat[][]{\label{doublespherea} \includegraphics[scale=0.25]{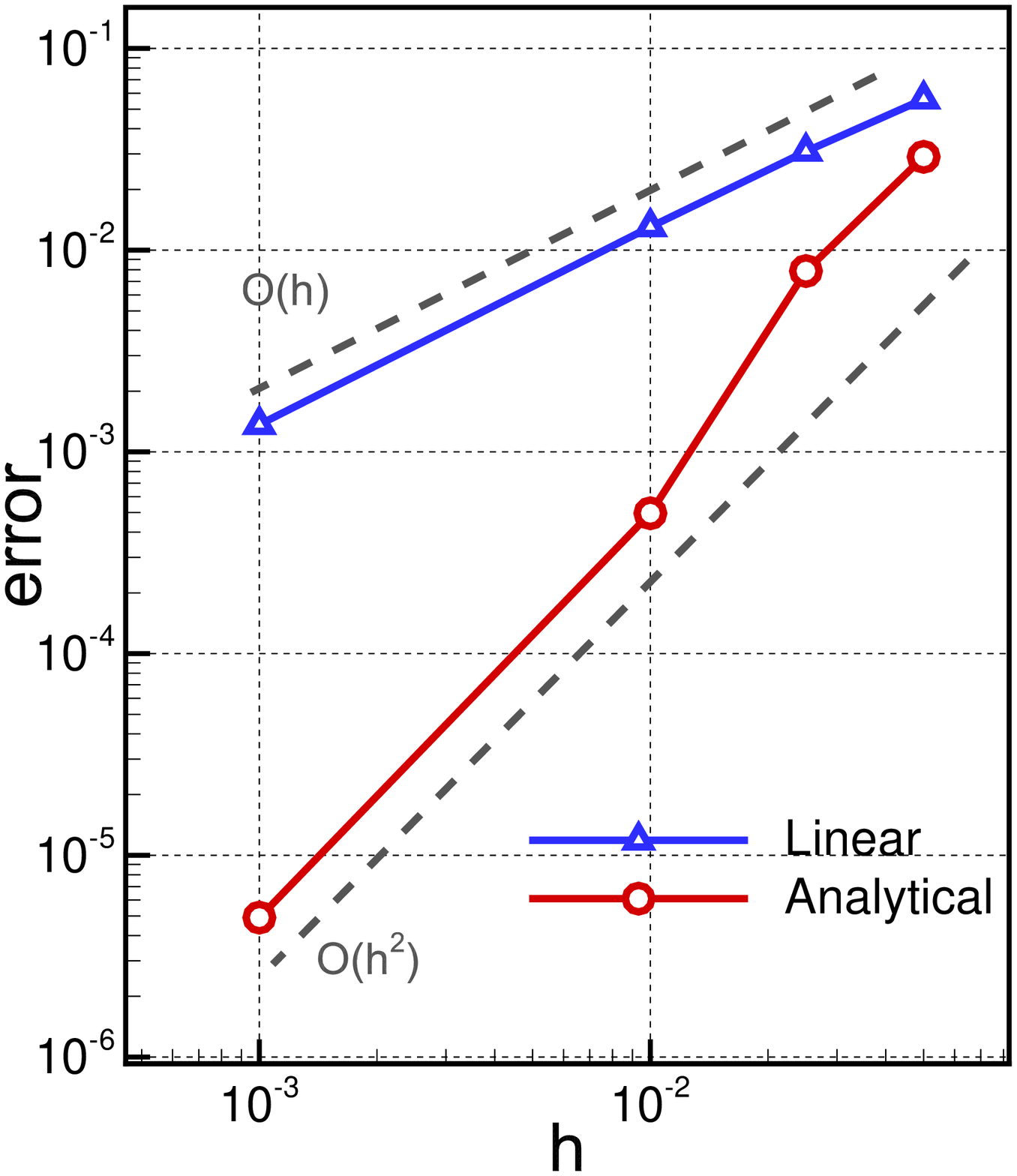}}
\subfloat[][]{\label{doublesphereb} \includegraphics[scale=0.25]{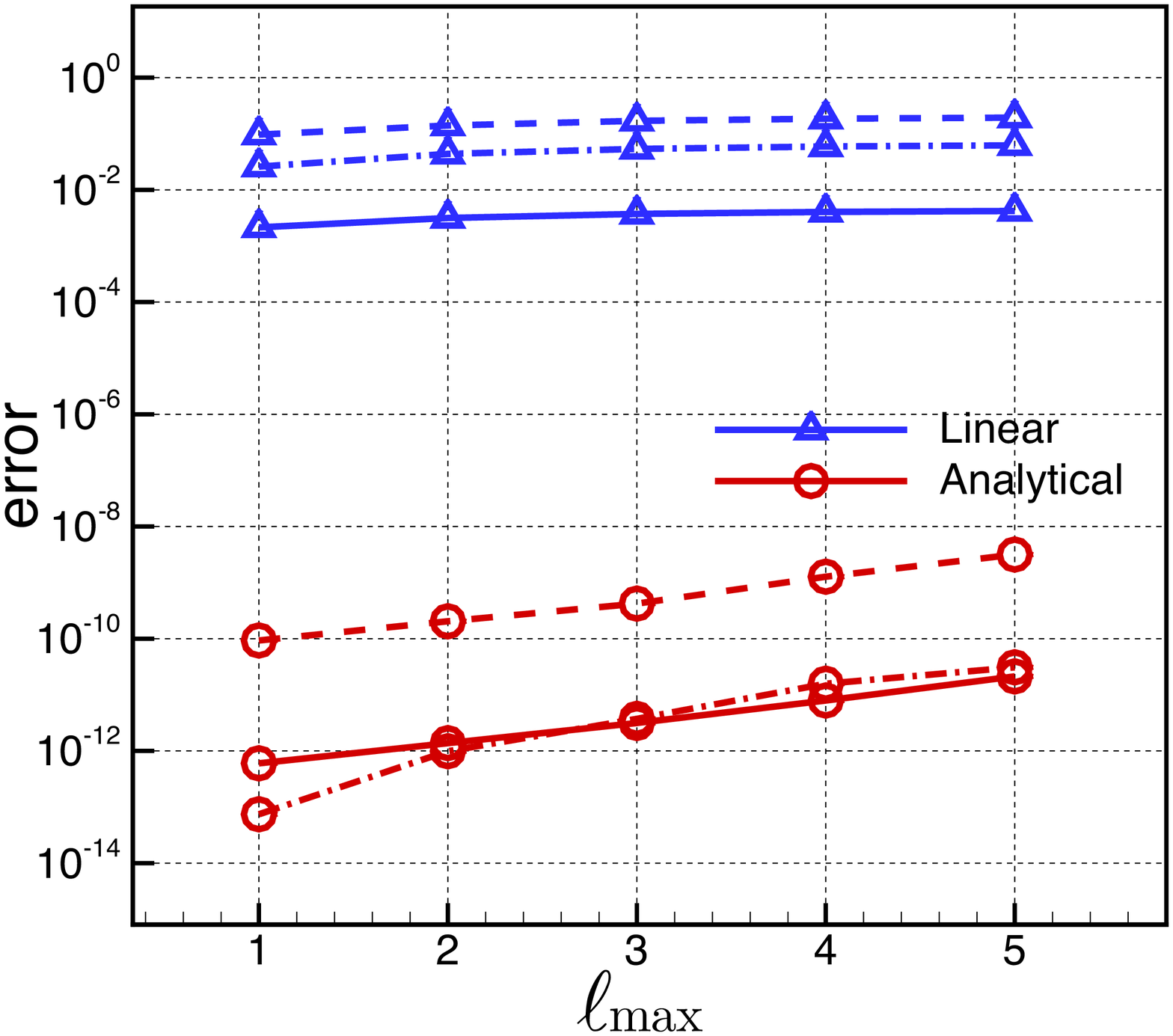}}
\caption{\label{doublesphere} Accuracy and consistency results for the 3D double-sphere case.}
\end{center}
\end{figure}

\begin{figure}
\begin{center}
\includegraphics[scale=0.25]{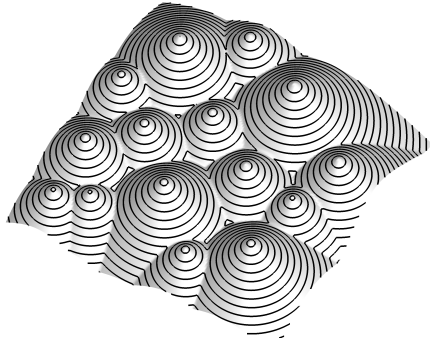}
\caption{\label{randomcontour} The contours of $15$ random generated circles.}
\end{center}
\end{figure}

In this section, we apply the analytical formulation to more complicated cases. Firstly we give two cases which have exact volumes to test both accuracy and consistency. The first one is refered to as ``double-circle'' case where two circles with an identical radius $r=0.25$ and different centers ($\textbf{x}_0 = (0.3, 0.5)$ and $\textbf{x}_1 = (0.7, 0.5)$) are merged. While the Zalesak disk with a radius of $0.4$, a notch width of $0.2$ and a notch height of $0.6$ is placed in at $(0.5, 0.5)$. As shown in Figs. \ref{geo_resulta} and \ref{geo_resultc}, the results converge with a 2nd-order rate in our method and a 1st-order rate in the liner method. Both these cases show distinguished consistency in the analytical formulation and large inconsistency in the linear method. A 3D double-sphere case with the same parameters in the double-circle case confirms the conclusion, as shown in Fig. \ref{doublesphere}. In Fig. \ref{randomcontour}, A case with $15$ randomly generated circles provides similar convergence and consistence results with the above cases.

\begin{figure}
\begin{center}
\includegraphics[scale=0.3]{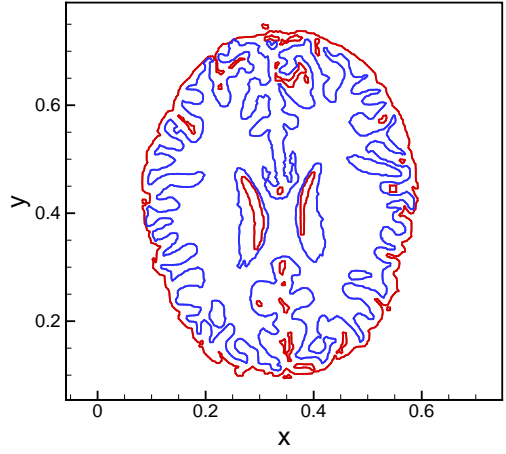}
\caption{\label{geobrain} The segmentation of a human brain generated by the method in Ref. \cite{li2011level}.}
\end{center}
\end{figure}

\begin{figure}
\begin{center}
\subfloat[][]{\label{braina} \includegraphics[scale=0.25]{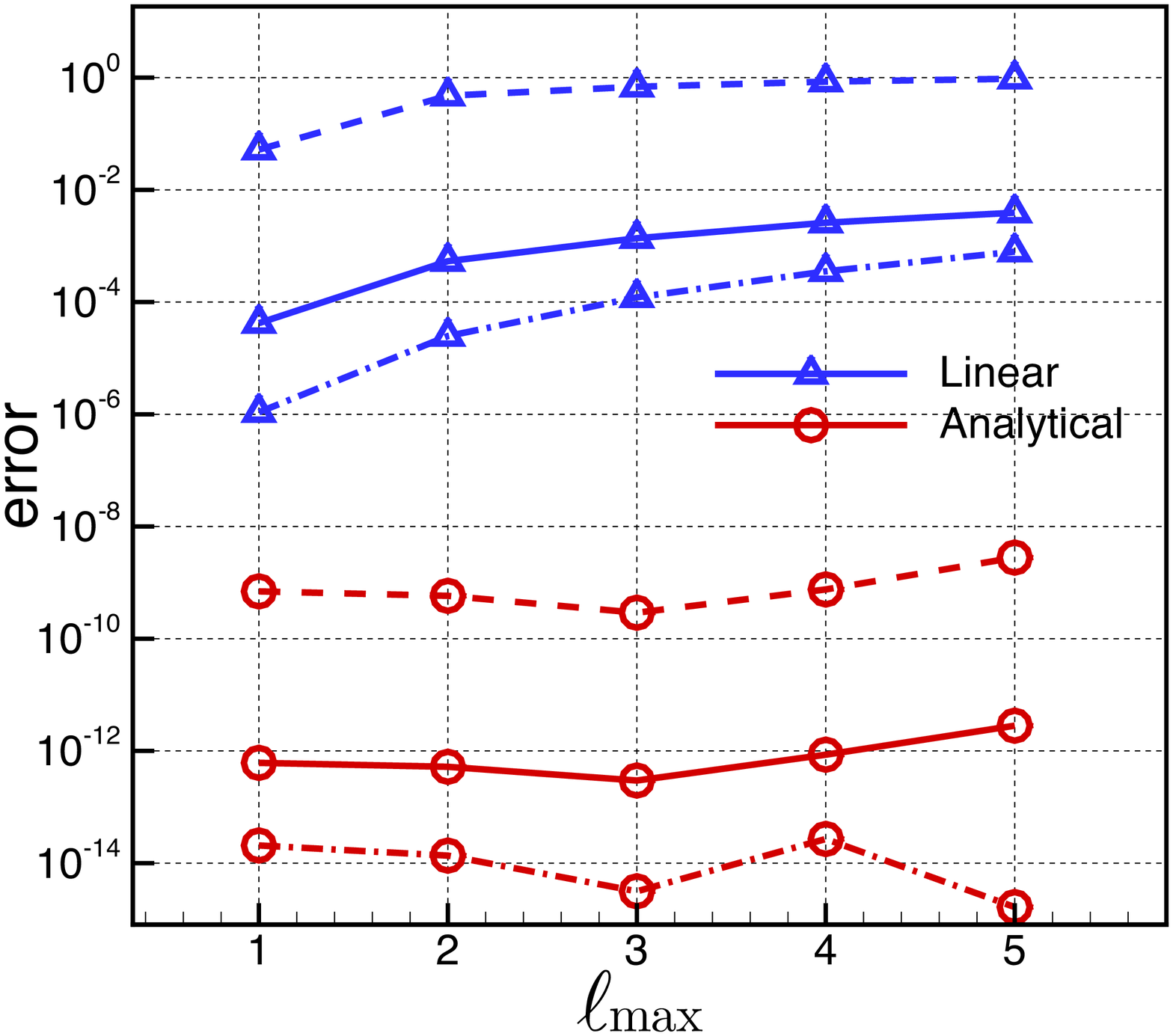}}
\subfloat[][]{\label{brainb} \includegraphics[scale=0.25]{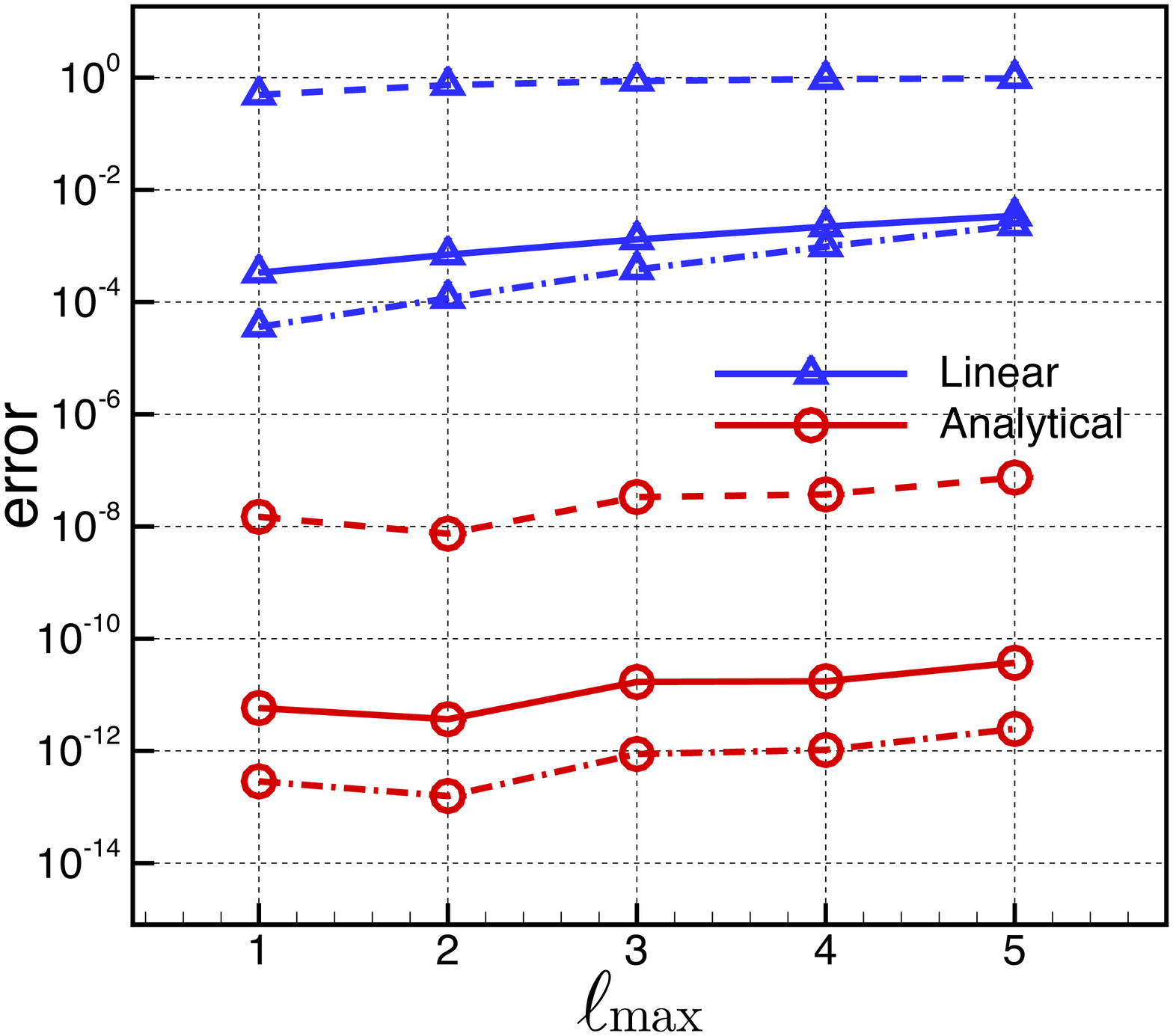}}
\caption{\label{brain} Consistency results of two different brain tissues for: (a) read line and (b) blue line in Fig. \ref{geobrain}.}
\end{center}
\end{figure}

Further more, a more practical case is tested to manifest the features of the present method. As shown in Fig. \ref{geobrain}, this case is the image segmentation of a human brain wihch has more complicated boundaries than the above cases. The level-set fields are generated by the code in Ref. \cite{li2011level} with the same MR image. The red line represents the air-tissue interface while the blue one stands for the interface of two different tissues. The consistency results are presented in Fig. \ref{brain}. Apparently, the errors in both areas are significantly reduced by the analytical formulation, indicating the consistency of our method in this brain image segmentation.
\section{Concluding remarks} \label{sec:con}
In this paper we have proposed a volume-estimation method for implicitly defined geometries. We have derived a general analytical formulation for 2D cases and three elementary 3D cases which can be employed for more complicated 3D cases. The analytical formulation achieves 2nd-order accuracy for arbitrary geometries. Based on this analytical formulation, the consistency issue during mesh refinement is imposed by the bilinear/trilinear interpolation. The accuracy and consistency are demonstrated by 2D and 3D examples. This consistent analytical volume-estimation method has the potential to be applied in many scientific and engineering fields, such as astronomy, computational fluid dynamics and medical imaging. Based on the idea developed in this paper, a consistent analytical volume-estimation formulation with high-order accuracy and its extension to arbitrary convex polyhedra is subject to further work.
\section*{Acknowledgment}
This work is supported by China Scholarship Council under No. 201306290030.
\section*{References}
\bibliographystyle{plain}
\bibliography{aipsamp}

\providecommand{\noopsort}[1]{}\providecommand{\singleletter}[1]{#1}%
\begin{thebibliography}{10}

\bibitem{besse20093}
S~Besse, O~Groussin, L~Jorda, P~Lamy, M~Kaasalainen, G~Gesquiere, E~Remy, and
  OSIRIS Team.
\newblock 3-dimensional reconstruction of asteroid 2867 {Steins}.
\newblock In {\em Lunar and Planetary Science Conference}, volume~40, page
  1545, 2009.

\bibitem{Cheong2007}
Benjamin Cheong, Mario F.~Rubin Raja~Muthupillai, and Scott~D. Flamm.
\newblock Normal values for renal length and volume as measured by magnetic
  resonance imaging.
\newblock {\em Clinical Journal of the American Society of Nephrology},
  2(1):38--45, 2007.

\bibitem{coakley2002prostate}
Fergus~V Coakley, John Kurhanewicz, Ying Lu, Kirk~D Jones, Mark~G Swanson,
  Silvia~D Chang, Peter~R Carroll, and Hedvig Hricak.
\newblock {Prostate Cancer Tumor Volume: Measurement with Endorectal MR and MR
  Spectroscopic Imaging}.
\newblock {\em Radiology}, 223(1):91--97, 2002.

\bibitem{diot2016interface}
Steven Diot and Marianne~M Fran{\c{c}}ois.
\newblock An interface reconstruction method based on an analytical formula for
  {3D} arbitrary convex cells.
\newblock {\em Journal of Computational Physics}, 305:63--74, 2016.

\bibitem{diot2014interface}
Steven Diot, Marianne~M Fran{\c{c}}ois, and Edward~D Dendy.
\newblock An interface reconstruction method based on analytical formulae for
  {2D} planar and axisymmetric arbitrary convex cells.
\newblock {\em Journal of Computational Physics}, 275:53--64, 2014.

\bibitem{Durso2014}
Timothy~A Durso, Jonathan Carnell, Thomas~T. Turk, and Gopal~N. Gupta.
\newblock {Three-Dimensional Reconstruction Volume: A Novel Method for Volume
  Measurement in Kidney Cancer}.
\newblock {\em The Journal of urology}, 28(6):745--750, 2014.

\bibitem{fujiwara2006rubble}
Akira Fujiwara, J~Kawaguchi, DK~Yeomans, M~Abe, T~Mukai, T~Okada, J~Saito,
  H~Yano, M~Yoshikawa, DJ~Scheeres, et~al.
\newblock {The rubble-pile asteroid Itokawa as observed by Hayabusa}.
\newblock {\em Science}, 312(5778):1330--1334, 2006.

\bibitem{Gong2012}
In~Hyuck Gong, Jinho Hwang, Don~Kyung Choi, Seung~Ryeol Lee, Young~Kwon Hong,
  Jae~Yup Hong, Dong~Soo Park, and Hwang~Gyun Jeon.
\newblock Relationship among total kidney volume, renal function and age.
\newblock {\em The Journal of urology}, 187(1):344--349, 2012.

\bibitem{gueyffier1999volume}
Denis Gueyffier, Jie Li, Ali Nadim, Ruben Scardovelli, and St{\'e}phane
  Zaleski.
\newblock Volume-of-fluid interface tracking with smoothed surface stress
  methods for three-dimensional flows.
\newblock {\em Journal of Computational Physics}, 152(2):423--456, 1999.

\bibitem{han2015scale}
L.H. Han, X.Y. Hu, and N.A. Adams.
\newblock Scale separation for multi-scale modeling of free-surface and
  two-phase flows with the conservative sharp interface method.
\newblock {\em Journal of Computational Physics}, 280:387--403, 2015.

\bibitem{hu2010multi}
XY~Hu, NA~Adams, M~Herrmann, and G~Iaccarino.
\newblock Multi-scale modeling of compressible multi-fluid flows with
  conservative interface method.
\newblock In {\em Proceedings of the Summer Program, Center for Turbulence
  Research}, 2010.

\bibitem{hu2006conservative}
X.Y. Hu, B.C. Khoo, N.A. Adams, and F.L. Huang.
\newblock A conservative interface method for compressible flows.
\newblock {\em Journal of Computational Physics}, 219(2):553--578, 2006.

\bibitem{joe1999brain}
Bonnie~N Joe, Melanie~B Fukui, Carolyn~Cidis Meltzer, Qing-shou Huang, Roger~S
  Day, Phil~J Greer, and Michael~E Bozik.
\newblock {Brain Tumor Volume Measurement: Comparison of Manual and
  Semiautomated Methods}.
\newblock {\em Radiology}, 212(3):811--816, 1999.

\bibitem{Jones1983}
T.~B. Jones, L.~R. Riddick, M.~D. Harpen, R.~L. Dubuisson, and D~Samuels.
\newblock Ultrasonographic determination of renal mass and renal volume.
\newblock {\em Journal of Ultrasound in Medicine}, 2(4):151--154, 1983.

\bibitem{lauer2012numerical}
E~Lauer, XY~Hu, Stefan Hickel, and Nikolaus~A. Adams.
\newblock Numerical modelling and investigation of symmetric and asymmetric
  cavitation bubble dynamics.
\newblock {\em Computers \& Fluids}, 69:1--19, 2012.

\bibitem{li2011level}
Chunming Li, Rui Huang, Zhaohua Ding, J~Chris Gatenby, Dimitris~N Metaxas, and
  John~C Gore.
\newblock {A level set method for image segmentation in the presence of
  intensity inhomogeneities with application to MRI}.
\newblock {\em IEEE Transactions on Image Processing}, 20(7):2007--2016, 2011.

\bibitem{Lorensen1987marchingcube}
William~E Lorensen and Harvey~E. Cline.
\newblock Marching cubes: A high resolution {3D} surface construction
  algorithm.
\newblock {\em ACM siggraph computer graphics}, 21(4), 1987.

\bibitem{mayr1996usefulness}
Nina~A Mayr, Vincent~A Magnotta, James~C Ehrhardt, James~A Wheeler, Joel~I
  Sorosky, B-Chen Wen, Charles~S Davis, Retta~E Pelsang, Barrie Anderson,
  J~Fred Doornbos, David~H Hussey, and William~T.C. Yuh.
\newblock Usefulness of tumor volumetry by magnetic resonance imaging in
  assessing response to radiation therapy in carcinoma of the uterine cervix.
\newblock {\em International Journal of Radiation Oncology* Biology* Physics},
  35(5):915--924, 1996.

\bibitem{Newman2006}
Timothy~S Newman and Hong Yi.
\newblock A survey of the marching cubes algorithm.
\newblock {\em Computers \& Graphics}, 30(5):854--879, 2006.

\bibitem{nystroem2002area}
Ingela Nystroem, Jayaram~K Udupa, George~J Grevera, and Bruce~E Hirsch.
\newblock Area of and volume enclosed by digital and triangulated surfaces.
\newblock In {\em Medical Imaging 2002}, pages 669--680. International Society
  for Optics and Photonics, 2002.

\bibitem{rypens2001fetal}
Fran{\c{c}}oise Rypens, Thierry Metens, Nathalie Rocourt, Pascale Sonigo,
  Francis Brunelle, Marie~Pierre Quere, Laurent Guibaud, Brigitte
  Maugey-Laulom, Chantal Durand, Fred~E Avni, et~al.
\newblock {Fetal Lung Volume: Estimation at MR Imaging
  $\rule[0.5ex]{1em}{0.5pt}$ Initial Results}.
\newblock {\em Radiology}, 219(1):236--241, 2001.

\bibitem{scardovelli2000analytical}
Ruben Scardovelli and St{\'e}phane Zaleski.
\newblock Analytical relations connecting linear interfaces and volume
  fractions in rectangular grids.
\newblock {\em Journal of Computational Physics}, 164(1):228--237, 2000.

\bibitem{so2011anti}
KK~So, XY~Hu, and NA~Adams.
\newblock Anti-diffusion method for interface steepening in two-phase
  incompressible flow.
\newblock {\em Journal of Computational Physics}, 230(13):5155--5177, 2011.

\bibitem{sussman1999adaptive}
Mark Sussman, Ann~S Almgren, John~B Bell, Phillip Colella, Louis~H Howell, and
  Michael~L Welcome.
\newblock An adaptive level set approach for incompressible two-phase flows.
\newblock {\em Journal of Computational Physics}, 148(1):81--124, 1999.

\bibitem{Sussman2000voflevelset}
Mark Sussman and Elbridge~Gerry Puckett.
\newblock A coupled level set and volume-of-fluid method for computing {3D} and
  axisymmetric incompressible two-phase flows.
\newblock {\em Journal of Computational Physics}, 162(2):301--337, 2000.

\end{thebibliography}

\begin{figure}
\begin{center}
\subfloat[][]{\label{randoma} \includegraphics[scale=0.25]{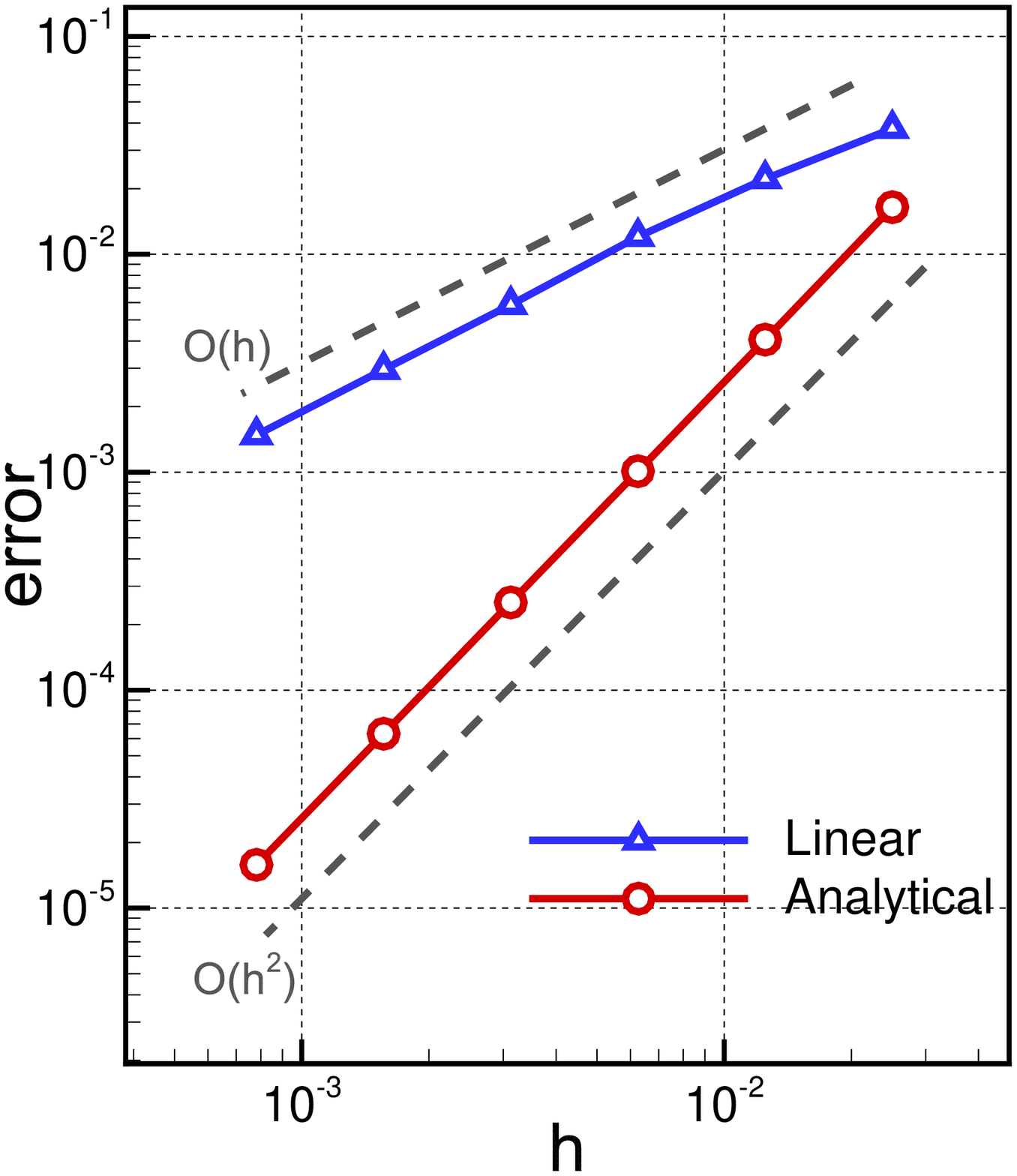}}
\subfloat[][]{\label{randomb} \includegraphics[scale=0.25]{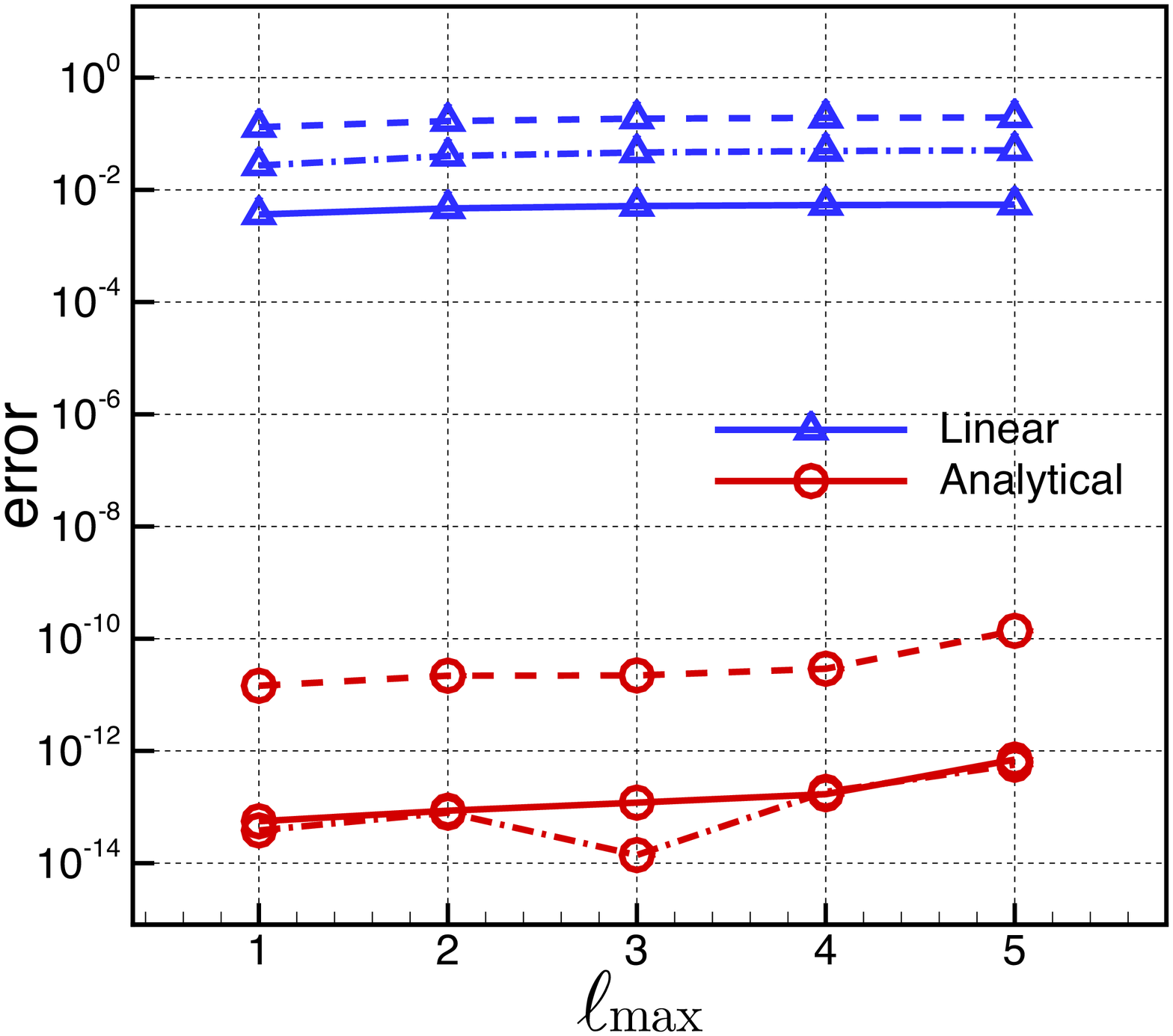}}
\caption{\label{random} Accuracy and consistency results for the random circles case.}
\end{center}
\end{figure}

\end{document}